\documentclass[journal]{new-aiaa}
\usepackage[utf8]{inputenc}
\usepackage{textcomp}

\usepackage{graphicx}
\usepackage{amsmath}
\usepackage[version=4]{mhchem}
\usepackage{siunitx}
\usepackage{longtable,tabularx}
\usepackage{subfig}
\usepackage{subcaption}
\usepackage{caption}
\usepackage{matlab-prettifier}
\usepackage{booktabs}
\usepackage{tabularx} 
\usepackage{multirow}
\usepackage{float}
\usepackage[inkscapelatex=false]{svg}
\usepackage[utf8]{inputenc}
\DeclareMathOperator*{\argmax}{\arg\!\max}

\newcolumntype{b}{X}

\title{Load Identification in Bistable Spacecraft Booms \\ via Parametric Data-Driven Modeling}

\author{Deven H. Mhadgut \footnote{Graduate Student, Department of Aerospace and Ocean Engineering}, Austin Phoenix\footnote{Director, Mission Systems Division, National Security Institute}, Serkan Gugercin\footnote{Professor, Department of Mathematics}, Samantha Parry Kenyon\footnote{Assistant Professor, Department of Aerospace and Ocean Engineering}}
\affil{Virginia Polytechnic Institute and State University, Blacksburg, VA 24060, USA}
\author{Jonathan Black\footnote{Flight Dynamics Thread Lead, Payload and Ground Systems Division}}
\affil{Northrop Grumman Corporation, Rolling Meadows, IL 60008 USA}
\author{Linus Balicki\footnote{Machine Learning Research Scientist}}
\affil{Novateur Research Solutions, Ashburn, VA 20147 USA}

\begin{document}

\maketitle

\begin{abstract}
Bistable tape spring booms are used on spacecraft for their ability to self-deploy using stored strain energy. However, their uncontrolled deployment can induce mechanical shocks that are variable as a function of material properties and temperature, and may damage sensitive satellite components and disrupt attitude control. Because traditional Finite Element Analysis (FEA) struggles to accurately capture this highly nonlinear behavior, we solve the inverse problem to estimate these loads from dynamic response measurements. 
Previous data-driven approaches using Vector Fitting required time-consuming retesting for every specific load level due to the boom's load-dependent dynamic behavior. To overcome this limitation, we introduce a parametric data-driven framework where a parametric transfer-function model of a composite tape spring boom is developed using force and velocity measurements. The parametric Adaptive Antoulas-Anderson algorithm (p-AAA) is used to construct a single parametric (multivariate) transfer function capable of capturing the nonlinear response of the boom to load amplitude. To evaluate the proposed framework, the boom is excited at its base at 15 distinct load levels using a single-axis reference input signal. Results demonstrate that the single parametric model outperformed the best discrete non-parametric case, reducing the total relative force estimation error for the reference signal by nearly 38\%. For experimental validation, the boom is subjected to sinusoidal, triangular and square signals. The cross validation results further supported this generalized performance. Collectively, these results show that the proposed parametric model accurately reconstructs input forces from velocity measurements alone, offering a solution for onboard diagnostics in future space missions.

\end{abstract}

\section{Introduction}\label{intro}

Bistable tape springs allow for compact, self-deploying structures in space applications; however, their rapid, uncontrolled deployment can produce unknown mechanical forces that may harm delicate components. These events are highly nonlinear, transient, impulsive and dependent on the operating environment. Measuring these deployment forces directly is often impractical due to challenges in placing a force sensor precisely at the point of force application and limitations in data transmission, particularly for small satellites. Furthermore, accurately characterizing these deployment forces during ground testing is further hindered by the difficulty in replicating the space environment. Previous studies have highlighted the risks associated with end-of-deployment shocks and complex behaviors like buckling and large rotations \cite{mallikarachchi_2011}, while others have shown that friction and viscoelastic relaxation from prolonged stowage can change deployment loads by as much as 30\% in around 3 hours at 60\textdegree C \cite{Mao2017DeploymentSystem,Kwok2013FoldingSprings,Brinkmeyer2016}. Similar effects have been observed due to temperature-dependence of the material properties with the natural frequency dropping with an increase in stowage temperature \cite{Adamcik2020, Mhadgut2023, Mhadgut2024, DENG2024, Soykasap2009}. Recent research has sought to capture these complex dynamics with greater fidelity. Tortorici et al. \cite{Tortorici2024} utilized a free-floating platform to experimentally validate how the rapid uncoiling of bistable booms induces significant attitude perturbations on small satellites. Despite these advancements in modeling, accurately predicting these deployment-induced forces and moments remains a challenge due to measurement difficulties and variability of the deployment system, though efforts have been made to mitigate them using damping materials and active control methods \cite{LIU201555,SHI20243993,KAMESH20121310}. 

Previously, Ma and Lin introduced inverse methods to estimate forces acting on cantilever beams, demonstrating how system models combined with measured vibration data could reconstruct unknown inputs \cite{Ma2000}. Subsequent studies \cite{Ma2003,Ma2004} expanded this methodology to nonlinear systems and emphasized the importance of model fidelity and robustness against measurement noise. Lin \cite{Lin2012} and Ma et al. \cite{Ma1998} further explored inverse dynamics in nonlinear contexts and impulsive loading scenarios, laying critical groundwork for analyzing transient, deployment-like events. Recently, advancements in inverse techniques have focused on enhancing stability and performance in ill-posed systems. Khoo et al. \cite{Khoo2014} evaluated the pseudo-inverse method across different sensor-to-DOF ratios, which is especially relevant in aerospace applications where instrumentation is constrained. Adaptive weighting inverse methods and system identification-based approaches further improve estimation reliability in noisy environments \cite{Liu2000}. Complementary work by Sethu et al. \cite{Sethu2015} and Steltzner and Kammer \cite{SteltznerKammer} compared various force reconstruction strategies, including frequency-domain filters and Frequency Response Function (FRF)-based techniques, offering insights into hybrid approaches that blend physical modeling with data-driven inference. Collectively, these studies provide a strong methodological basis for applying inverse force estimation techniques to deployment events in space structures, where direct force measurement remains infeasible. 

Previous studies have often relied on accurate physical models using Finite Element Analysis (FEA) \cite{mallikarachchi2014deployment, Shore2022}. Our prior research demonstrated that these systems are very sensitive to geometric uncertainty, which we mitigated by using 3D scanning techniques to reduce geometric errors in our physical models \cite{Mhadgut2025}. However, even with precise geometry, inaccuracies persist in higher modes due to physical uncertainties such as root boundary conditions and loading and material nonlinearities. Furthermore, the structural dynamics of bistable tape springs are highly sensitive to input load amplitude. As the load and vibration amplitude increases, the system exhibits severe nonlinearities causing the underlying frequency response functions to shift significantly \cite{shabeer2025}. Previous studies frequently assumed linear system behavior or focused on idealized excitations such as single-frequency sinusoids or impulses. These limitations reduce their applicability to real-world scenarios like the deployment of tape spring booms, which involve complex transient loads. Additionally, many approaches have not been validated experimentally across a wide range of excitation types or lacked integration with practical sensing strategies suitable for lightweight space structures. Our previous work \cite{Mhadgut2026} highlighted the load dependent nonlinearities by employing a data-driven rational approximation methodology that estimated the input force from the dynamic response of the boom. It was observed that this rational approximation modeling was limited by its inability to generalize for different load levels and retesting was needed to define frequency response (transfer) functions at each distinct load level. To overcome this limitation, this study introduces a parametric modeling approach using the parametric Adaptive Antoulas-Anderson (p-AAA) algorithm \cite{Rodriguez2023}. This approach leads to a single, unified, yet parametric,  
transfer function model capable of accurately capturing the system's load-dependent nonlinear behavior at any point within the tested load range. 

Several shaker tests are conducted on a composite tape-spring boom to obtain FRFs with a broadband periodic chirp reference signal as the input and the measured tip velocity from a laser vibrometer as the output. The FRF data is then used to construct a single parametric transfer-function model of the boom system which can be evaluated at various input load levels. For experimental validation of the model, we conduct more shaker tests where the boom is excited using various signal types, including periodic chirp, sinusoidal, triangle and square. Input forces are estimated from the measured responses and compared to the known input excitations for validation. This process enables a systematic evaluation of the estimation technique’s accuracy under controlled conditions while replicating the dynamics relevant to real deployment scenarios. Essentially, this effort focuses strictly on the initial development and experimental validation of a high-accuracy, force-dependent parametric model. The subsequent step of using this framework to calculate actual deployment input loads is beyond the scope of this paper and will be addressed in future work.

The paper is organized as follows: Section \ref{boommat} presents the deployer design and boom material properties, Section \ref{exp} describes the experimental setup and Section \ref{colmet} outlines the various estimation and data-processing techniques used in this work. Finally, the results are presented and discussed in Section \ref{parametric}, and the paper concludes in Section \ref{conclusion}.

\section{Boom Deployer Design and Material Properties of the Boom}\label{boommat}

\begin{figure}[h]
\centering
\begin{tabular}{cc}
     \subfloat[Boom deployer render with boom stowed inside \cite{engebretson2022hybrid} \label{dep}]{\includegraphics[width=0.46\textwidth]{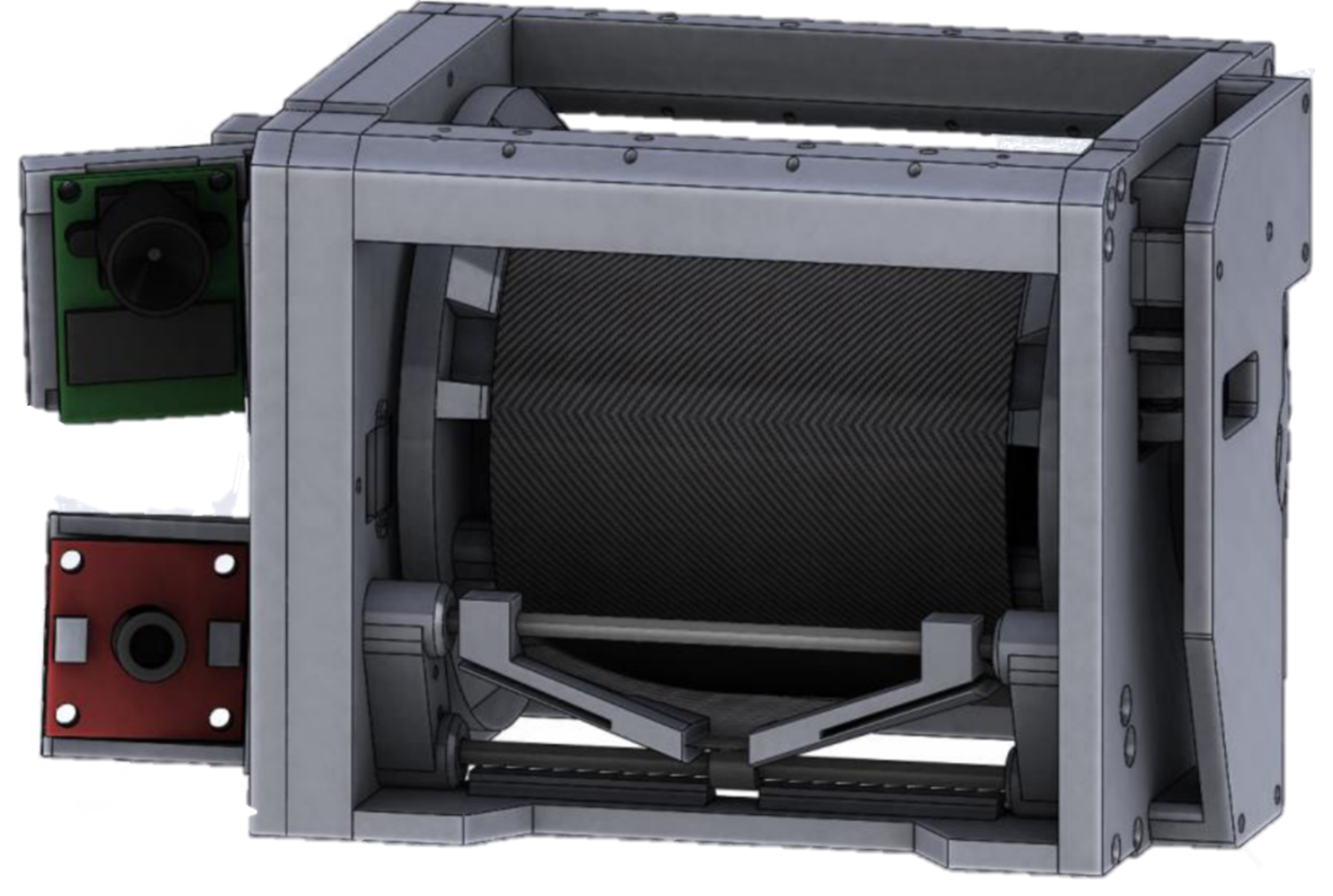}}
      & 
      \subfloat[Boom cross-section showing three plies and embedded copper cables \label{boomcs}]{\includegraphics[width=0.39\textwidth]{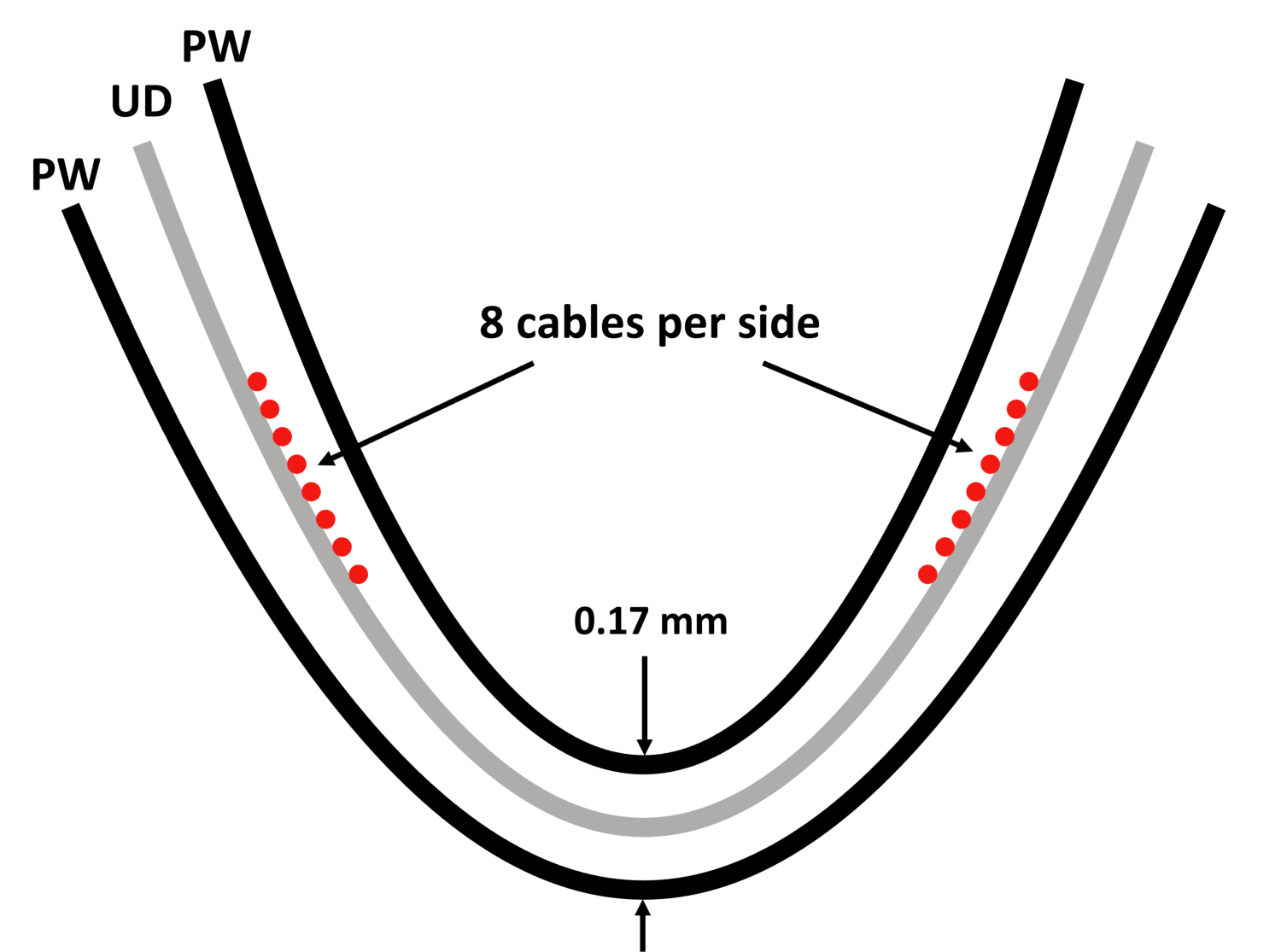}}
\end{tabular}
\caption{Boom deployer and composite layup}
\end{figure}

The boom used in the experiments was fabricated at NASA Langley and is identical to the one scheduled for deployment on Ut ProSat-1 (UPS-1), a 3U CubeSat developed by students at Virginia Tech \cite{Engebretson2022UtBooms}. The primary objective of the UPS-1 mission is to characterize the boom's performance throughout the mission by deploying a 4-foot boom and subsequently retracting it using a stepper motor while capturing the satellite’s dynamic response, and repeating the sequence multiple times along the orbit. The mission is designed to evaluate the boom’s performance under various operational conditions throughout its lifetime, serving as a verification test for the deployment system. The structural response is recorded via two Inertial Measurement Units (IMUs), one positioned at the boom tip \cite{Yao2024} and the other at the root within the CubeSat chassis. The tip-mounted IMU is integrated into a flexible circuit connected using copper cables that run through the boom (shown in Fig. \ref{boomcs}). Deployment is facilitated by a custom deployer (shown in Fig. \ref{dep}) incorporating a camera shutter design, wherein a servo motor restrains the spool to prevent premature deployment and releases it when needed \cite{engebretson2022hybrid}. 

The boom specimen is 4 ft long ($\sim$122 cm), 70 mm wide, and has a parabolic cross-section with a measured thickness of 0.17 mm. It features a [45PW/0UD/45PW] composite layup designed for reliable self-deployment after long-term stowage. The laminate consists of two outer ±45° carbon fiber/epoxy plain weave (M30S/PMT-F7) plies and one inner 0° unidirectional (MR60H/PMT-F7) ply. Key material properties of the composite are listed in Table \ref{tab:1}.

\begin{table*}[h]
\caption[Table]{Material properties of the thin-ply composite \cite{Lee2024}}\label{tab:1}
\centering{
\resizebox{\columnwidth}{!}
{
\begin{tabular}{c c c c c c c}
\toprule\toprule
Ply Material & Fiber/Resin & $E_1$ [GPa] & $E_2$ [GPa] & $\nu_{12}$ &  $G_{12}$ [GPa] & Thickness $t$ [$\mu$m]\\
\midrule
Unidirectional Carbon Fiber   & MR60H/PMT-F7     &   144.1 & 5.2 & 0.335 & 2.8 & 40.0 \\
Plain Weave Carbon Fiber  &  M30S/PMT-F7    &  89.0 & 89.0 & 0.035 & 4.2 & 58.2 \\
\bottomrule\bottomrule
\end{tabular}
}
}
\end{table*}

\section{Experimental Setup}\label{exp} 
\vspace{3 mm}

\begin{figure}[h!]
\centering
\includegraphics[width=\textwidth]{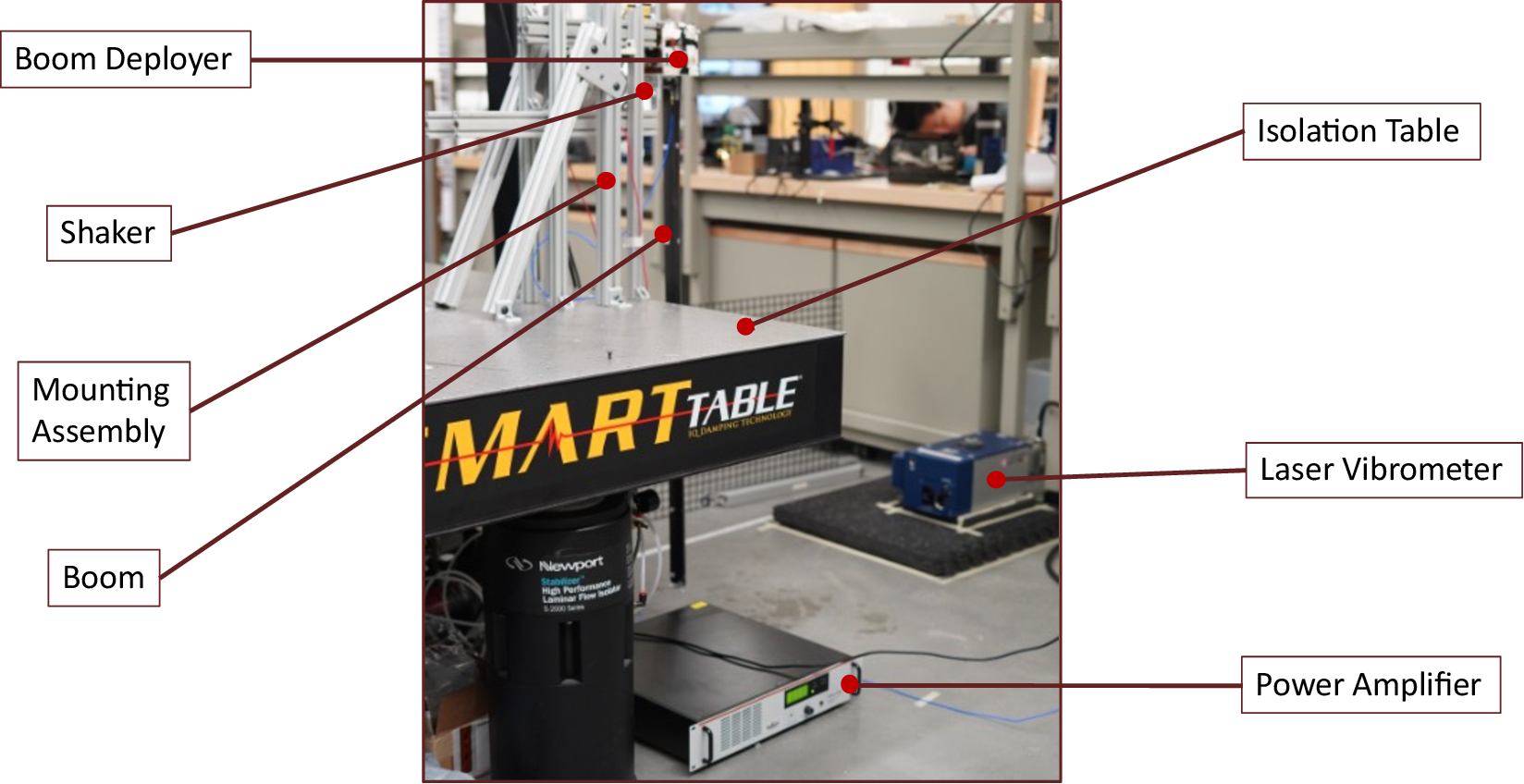}
\caption{Experimental setup}
\label{exp_setup_1}
\end{figure}

\begin{figure}[htbp]
\centering
\begin{tabular}{c c}

     \subfloat[Side View of the Boom Deployer attached to an Al Frame]{\includegraphics[width=0.32\textwidth]{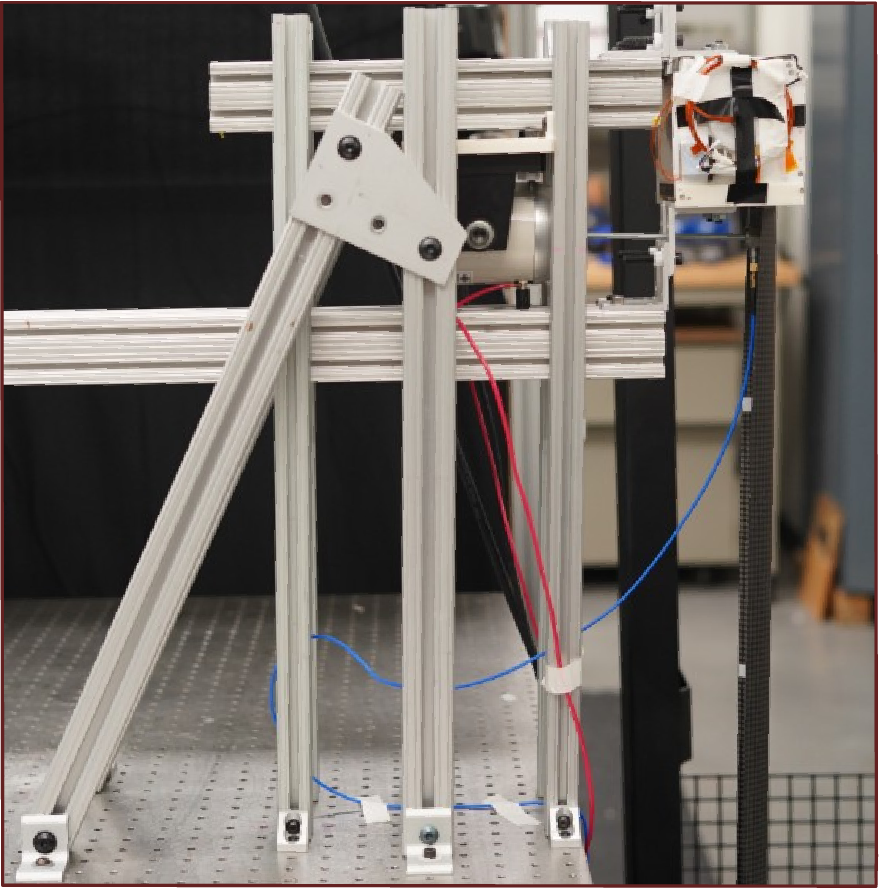}\label{sideview}} &
       
     \subfloat[Isometric View of the Boom Deployer attached to an Al Frame]{\includegraphics[width=0.3225\textwidth]{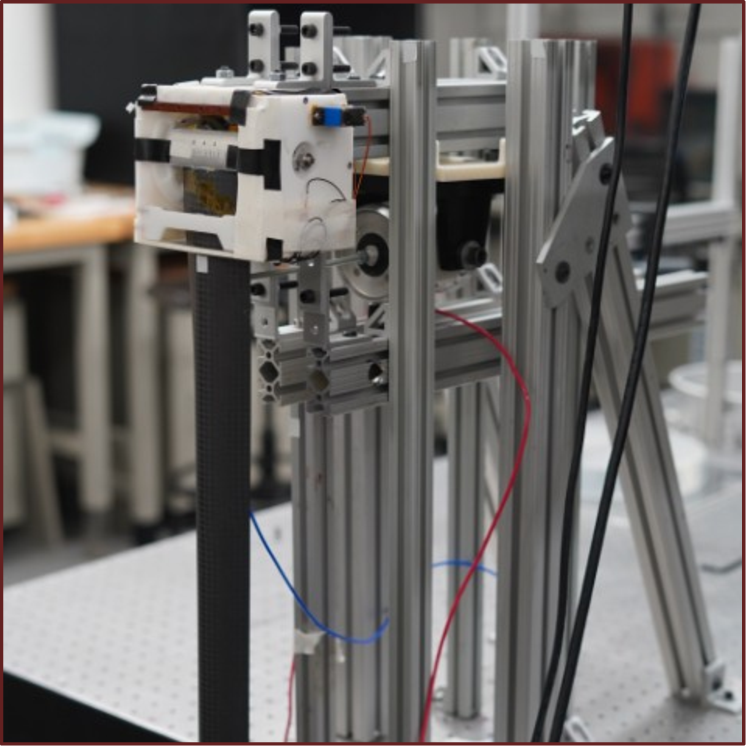}\label{isoview}}\\

     \subfloat[Shaker Input]{\includegraphics[width=0.32\textwidth]{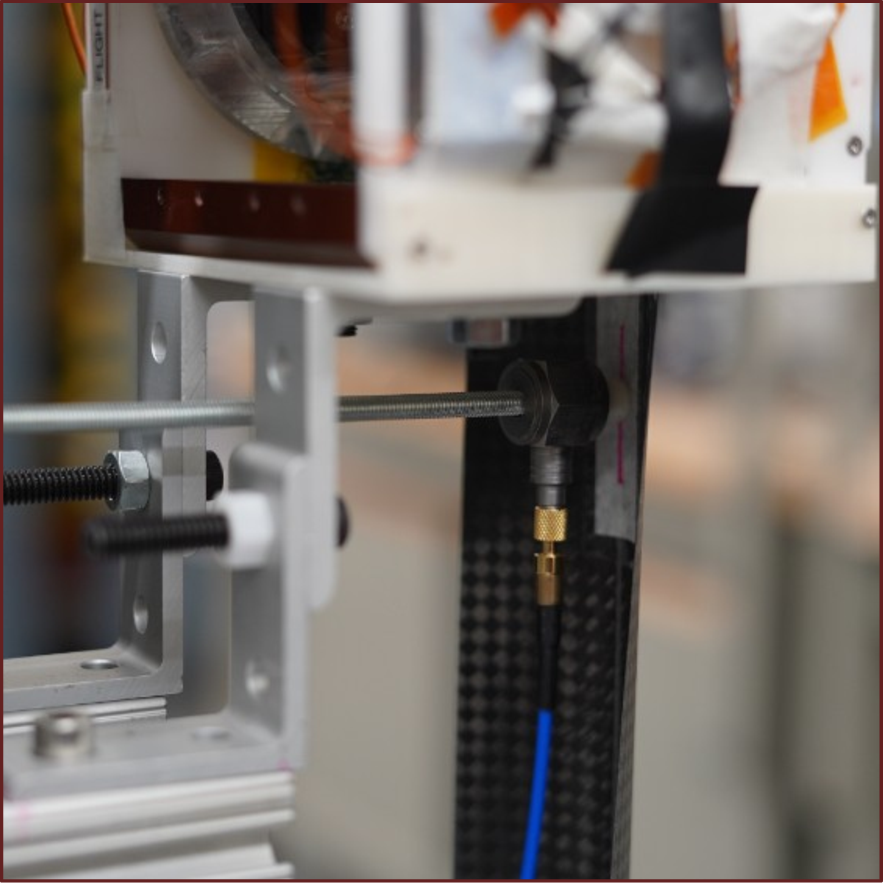}\label{input}} &

     \subfloat[Velocity Output near the tip of the boom]{\label{mode4}\includegraphics[width=0.323\textwidth]{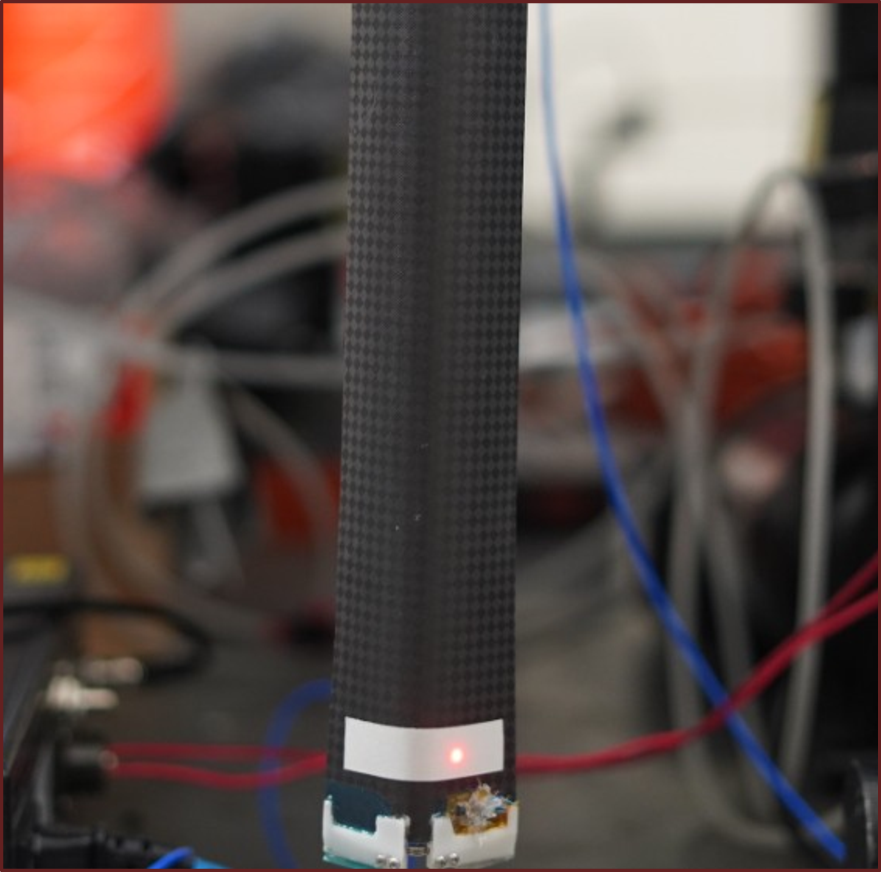}\label{output}}
  
\end{tabular}
\caption{Close-up of test setup}
\label{closeup}
\end{figure}

The experimental setup shown in Fig. \ref{exp_setup_1}, consists of the boom deployer fixed to a vibration isolation table using an aluminum extrusion frame. A 3D-printed clamp is used to secure the deployer to the test frame. The test frame was reinforced with additional supports from the back to minimize interference with the test results by having natural frequencies higher than the test frequency range. An electrodynamic mini-shaker (shown in Figs. \ref{sideview} and \ref{isoview}), rigidly attached to the test frame using a 3D-printed PLA fixture, was used to impart loads to the boom. The shaker stinger was aligned with the boom’s central axis near its base and connected via a load cell to the boom surface to measure the force imparted by the shaker at that location. A custom 3D-printed interface cap was affixed to the load cell to precisely conform to the boom's curvature to ensure accurate and uniform load transmission. Then, the load cell interface cap was secured to the boom surface using some super glue and double-sided tape to ensure a stable coupling. A close-up of the shaker input assembly is shown in Fig. \ref{input}. 

\begin{figure}[htpb]
\centering
\includegraphics[width=0.7\textwidth]{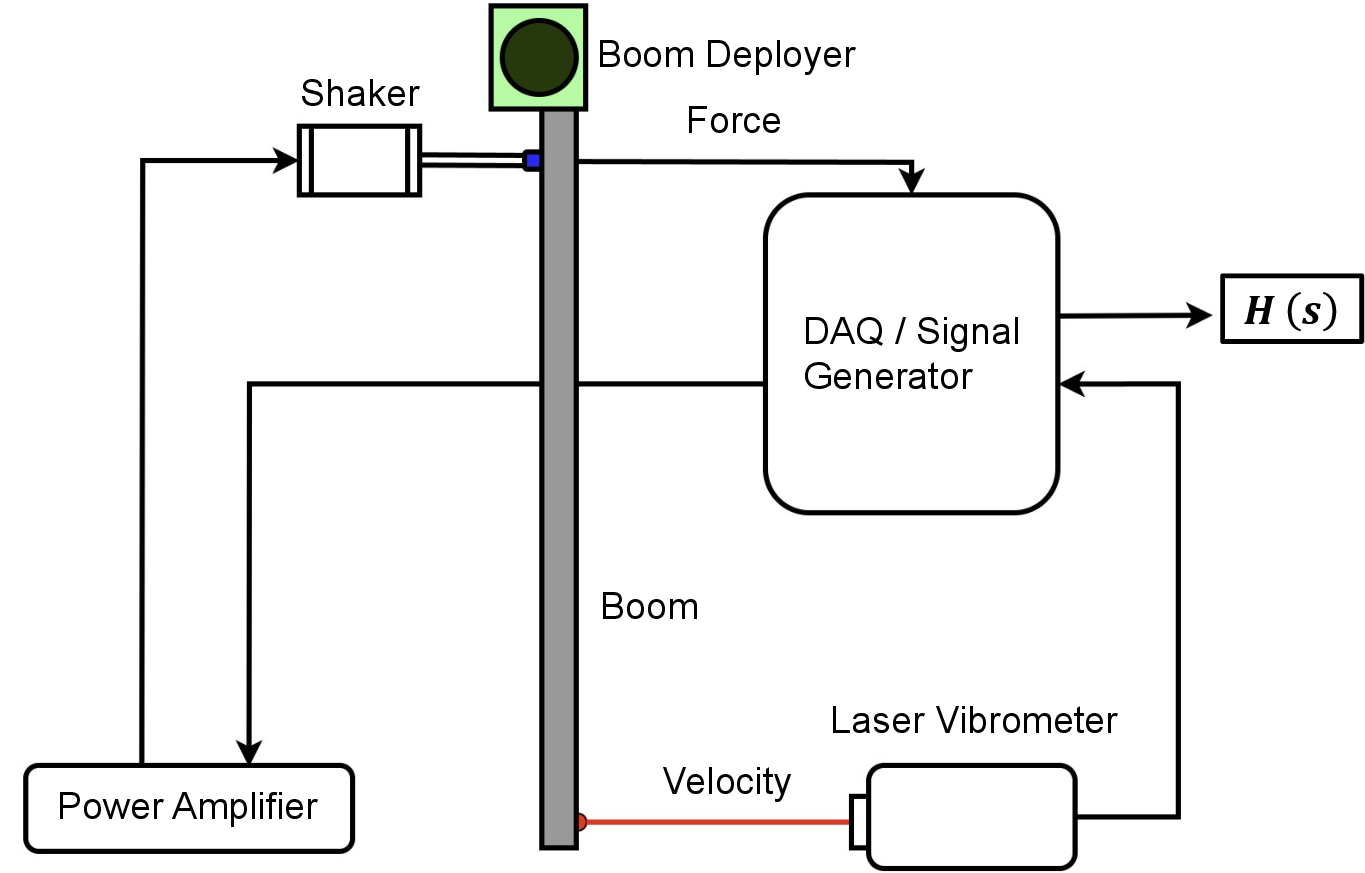}
\caption{Test setup schematic showing data flow}
\label{exp_schematic}
\end{figure}

A simplified test schematic showing the flow of data has been presented in Fig. \ref{exp_schematic}. The internal waveform generator of the Polytec Data Acquisition (DAQ) System generates a periodic chirp reference signal from 0.01 Hz to 1000 Hz and sends it to the shaker via a power amplifier. The shaker imparts the input force to the boom via a load cell (shown in blue). A Polytec PSV-400 laser vibrometer measures the velocity response at a point 15 mm from the tip of the boom (shown in Fig. \ref{output}) near the location where the IMU is attached which will be used for on-orbit acceleration measurements on UPS-1.  The Polytec DAQ system then analyzes the input and output data and generates frequency response function (FRF) data. Two types of transfer functions are computed: TF1, defined between the voltage input and the force output for the shaker; and TF2, between the force output from the shaker and the velocity measured at the tip of the boom using the vibrometer. In this paper, we only focus on TF2 because the primary objective is to characterize the structural dynamics of the boom rather than the input-output characteristics of the shaker system. Data is collected at a sampling rate of 256 Hz for 128 seconds yielding frequency response function samples at multiple frequencies, $s_j = \imath\omega_j$ for $j=1,2,\ldots,N_s$, where $N_s=12800$. To minimize the effect of noise on the results, complex averaging is used to obtain the $H_1$ estimator of the FRF given by
\begin{equation}\label{eq:1}
{H}_1(\imath\omega)=\frac{\mathcal{G}_{vf}(\imath\omega)}{\mathcal{G}_{ff}(\imath\omega)},
\end{equation}
where $\mathcal{G}_{\upsilon f}(\imath\omega)=E\left[V(\imath\omega)F(\imath\omega)^{*}\right]$ is the cross-spectral density and $\mathcal{G}_{ff}(\imath\omega)=E\left[F(\imath\omega)F(\imath\omega)^{*}\right]$ is the auto-spectral density. Here $E[\cdot]$ is the estimator function and $\{\cdot\}^{*}$ is the complex-conjugate of a matrix. This work utilizes the $H_1$ estimates of the FRFs to develop single-input single-output (SISO) data-driven models of the boom. For simplicity, the subscript is dropped and the experimental FRFs are denoted by $H (\imath\omega)$. Therefore, the experimental FRF data are: 
\begin{equation}\label{eq:2}
H\left(\imath\omega_j\right) \in \mathbb{C}^{1 \times 1} \text { for } j=1,2, \ldots, 12800.
\end{equation}

\section{Data Collection and Processing Methodology}\label{colmet}
\begin{figure}[htpb]
\centering
\includegraphics[width=\textwidth]{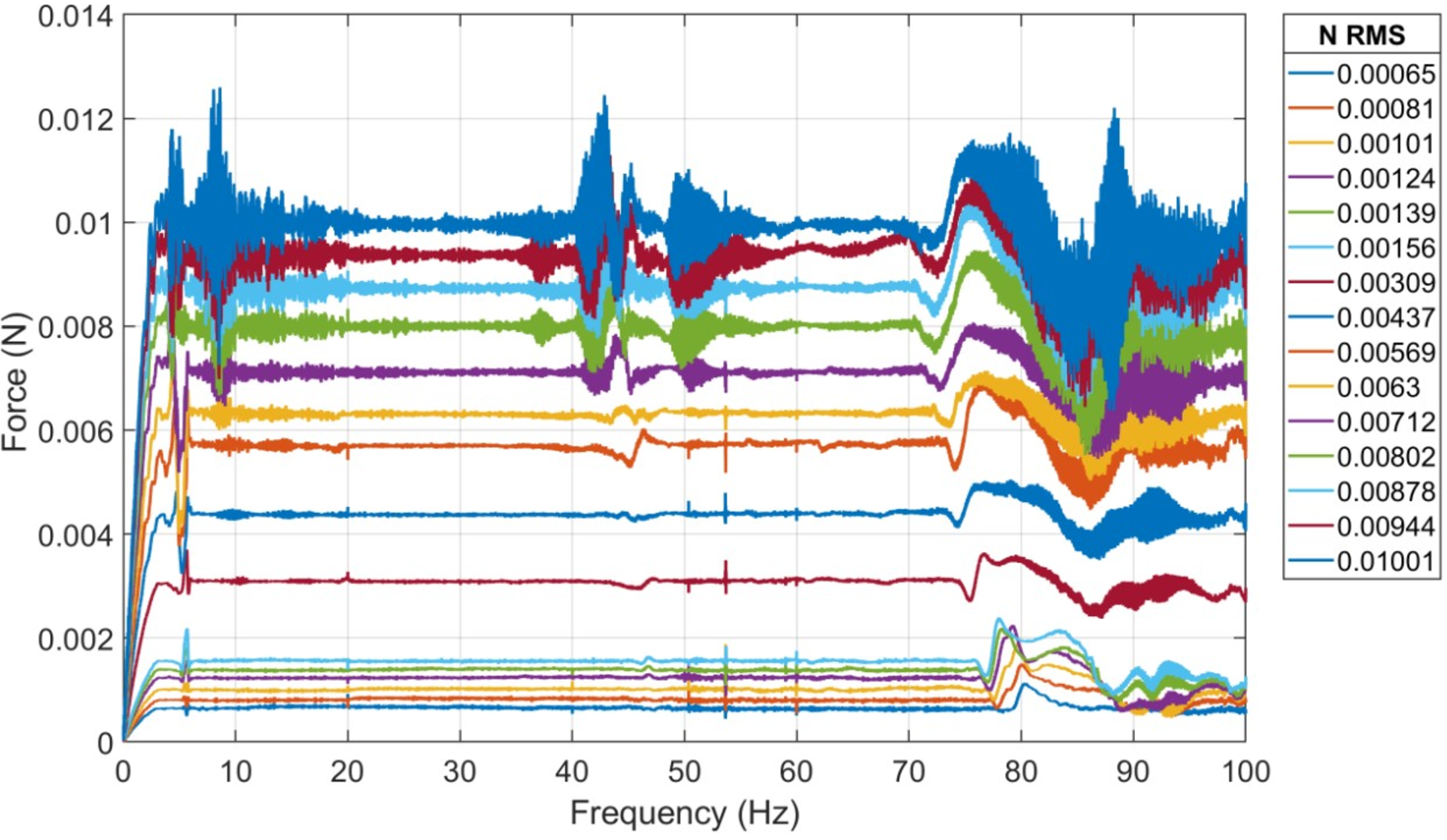}
\caption{Amplitude-corrected force input (periodic chirp signal) for all load cases in the frequency domain}
\label{ampcorrall}
\end{figure}

To assess the effect of input load, 15 load levels in the [0.00065, 0.01] N RMS range (in the frequency domain) were applied based on the shaker input voltages : [0.05, 0.06, 0.07, 0.08, 0.09, 0.1, 0.2, 0.3, 0.4, 0.5, 0.6, 0.7, 0.8, 0.9, 1.0] V. This range of input voltages is considered sufficient for eventually creating a parametrically trained model for predicting deployment loads. The methodology began with the collection of frequency response data from the boom system under test. The boom was excited using periodic chirp signals for a duration of 128 seconds. The input force spectra for all the tests are shown in Fig. \ref{ampcorrall}. The increased noise in the force input beyond 70 Hz is driven by dynamic interactions between the the boom and the deployer, where the nonlinearities of the boom become prominent. Retaining the full 0 to 100 Hz bandwidth is essential for two reasons. First, the data-driven model will be applied to an actual deployment, which is a highly nonlinear and transient event that excites a broad frequency spectrum. Second, this wider bandwidth would accurately capture the dynamic interactions between the boom and the deployer.

\begin{figure}[htbp]
\centering
\includegraphics[width=\textwidth]{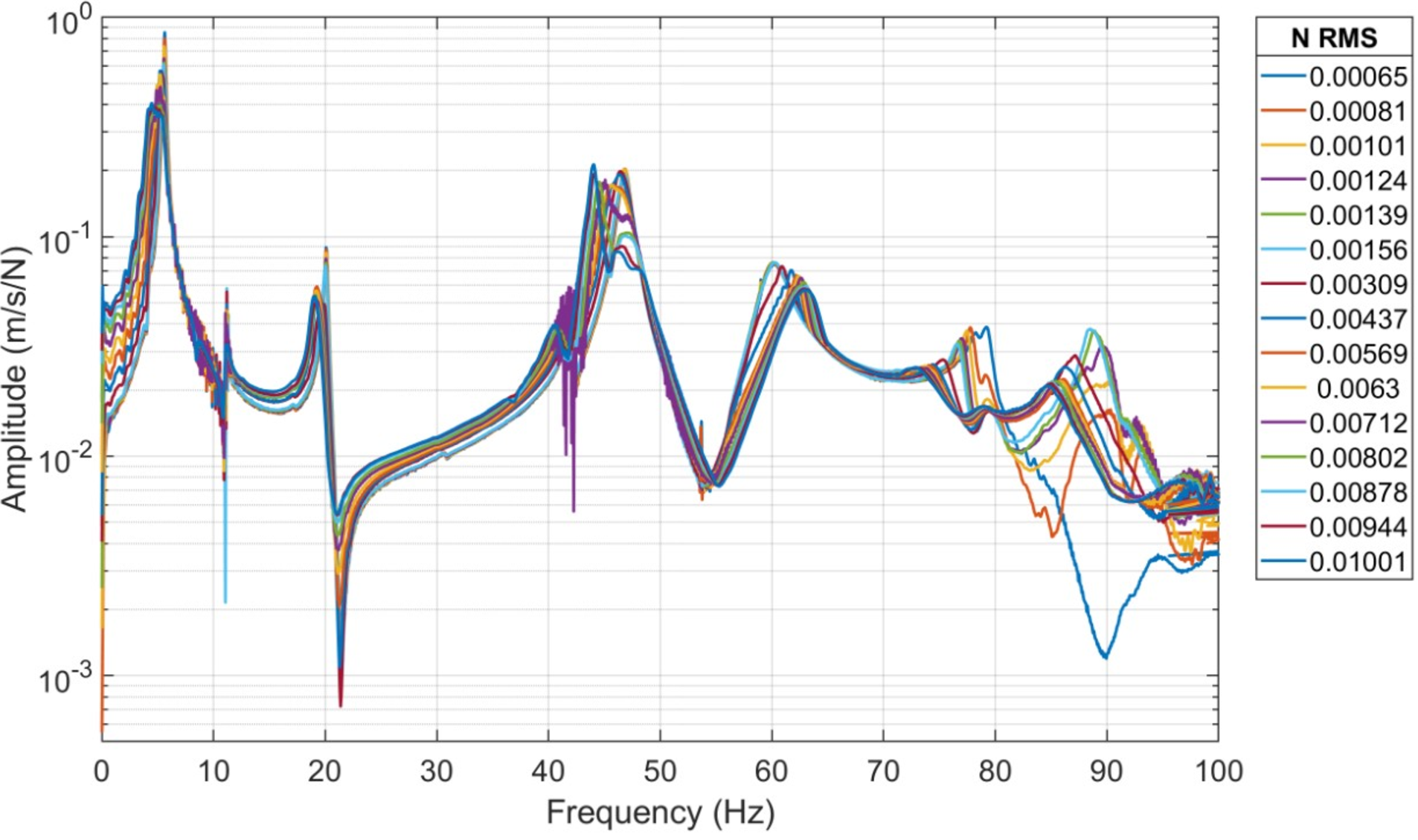}
\caption{Velocity/Force Transfer Functions }
\label{TF_all}
\end{figure}

Maintaining a flat unbiased input force profile in the frequency domain is inherently difficult for the shaker. This is because the mechanical impedance of the boom fluctuates drastically at its resonances and anti-resonances, which constantly changes how much force it draws from the shaker. However, the force spectra had relatively constant amplitude in the 3 to 70 Hz range. Therefore, the Root Mean Square (RMS) value of the force within this $3  - 70 \text{ Hz}$ frequency range is selected as the effective ``parameter", $p$, to provide a single-value representation for the different load levels. The legend in Figs. \ref{ampcorrall} and \ref{TF_all} shows the parameter values for the tests.

The force-velocity Frequency Response Functions (FRFs), presented in Fig. \ref{TF_all}, provide a comprehensive characterization of the system's dynamic behavior and serve as the basis for model identification. Inspection of these responses shows load-dependent shifts in both the natural frequency and damping ratio of the first bending mode. Furthermore, as the excitation force is applied normal to the boom surface, the first bending mode dominates the response and exhibits the highest amplitude in the transfer function data.

\subsection{Frequency Domain Transfer Function Inversion}\label{finv} 
\vspace{2 mm}
The dynamic relationship between an input force $f({t})$ and the resulting system velocity ${v}({t})$ can be represented in the time domain using the convolution integral:
\begin{equation}\label{Eq:3}
{v}\left(t\right)=\intop_0^th(t-\tau)f\left(\tau\right)d\tau,
\end{equation}
where $h(t)$ is the impulse response function of the system. Applying Fourier transform converts this convolution into a multiplication in the frequency domain:
\begin{equation}\label{eq:4}
{V}\left(\imath\omega\right)=[{H}\left(\imath\omega\right)]{F}(\imath\omega),
\end{equation}
where ${H}$ is the frequency response function (FRF). The FRF can be obtained experimentally using a known reference input force. Direct inversion of the experimental frequency response matrix $H$ in Eq.~\eqref{eq:4} amplifies noise, especially at anti-resonances where the magnitude approaches zero. To prevent this, the raw samples $H$ are first smoothed and a rational function is then fitted to these samples to create a continuous analytical model, $\tilde{H}$ for reliable force estimation. This fitting is done using Vector Fitting  and p-AAA approximation for the non-parametric and parametric studies, respectively. Then, using $\tilde{H}$ in place of the raw data $H$, Eq.~\eqref{eq:4} is solved to obtain the force signal in the frequency domain. The time history of the input force can then be obtained using the inverse Fourier transform, i.e.~
\begin{equation}\label{Eq:5}
    f(t) = \mathcal{F}^{-1} \{ {F}(\imath\omega)\}.
\end{equation}

\subsection{Non-Parametric Analysis and the Need for Parametric Modeling}

Following data acquisition, the frequency response data
$\{H(\imath\omega_j)\}$ in Eq.~\eqref{eq:2}, is subjected to \textit{Vector Fitting} \cite{Gustavsen1999}, an established technique to construct a rational approximant of the frequency-domain data.   
Data-driven rational approximation from measured data has been a heavily studied topic; see, e.g., \cite{aca90,sk1963,mayo2007framework,karachalios2020loewner,AntBG20,hokanson2017projected,Drmac-Gugercin-Beattie:VF-2014-SISC,berljafa2017rkfit,Nakatsukasa2018,gosea2020algorithms,ackermann2025frequency,Gosea2024Structured,benner2021model1} and the references therein. Following our previous work \cite{Mhadgut2026}, we employ the vector-fitting method  for the non-parametric rational approximation step; however, as stated earlier, the central contribution of the present work is the parametric extension via the p-AAA algorithm, which overcomes the load-level limitation and eliminates the need for retesting at each distinct load level. 

\begin{figure}[ht]
\centering
\begin{tabular}{c c}

     \subfloat[Load Case 1 : 0.0006 N $  \mathrm{RMS}$]{\includegraphics[width=0.38\textwidth]{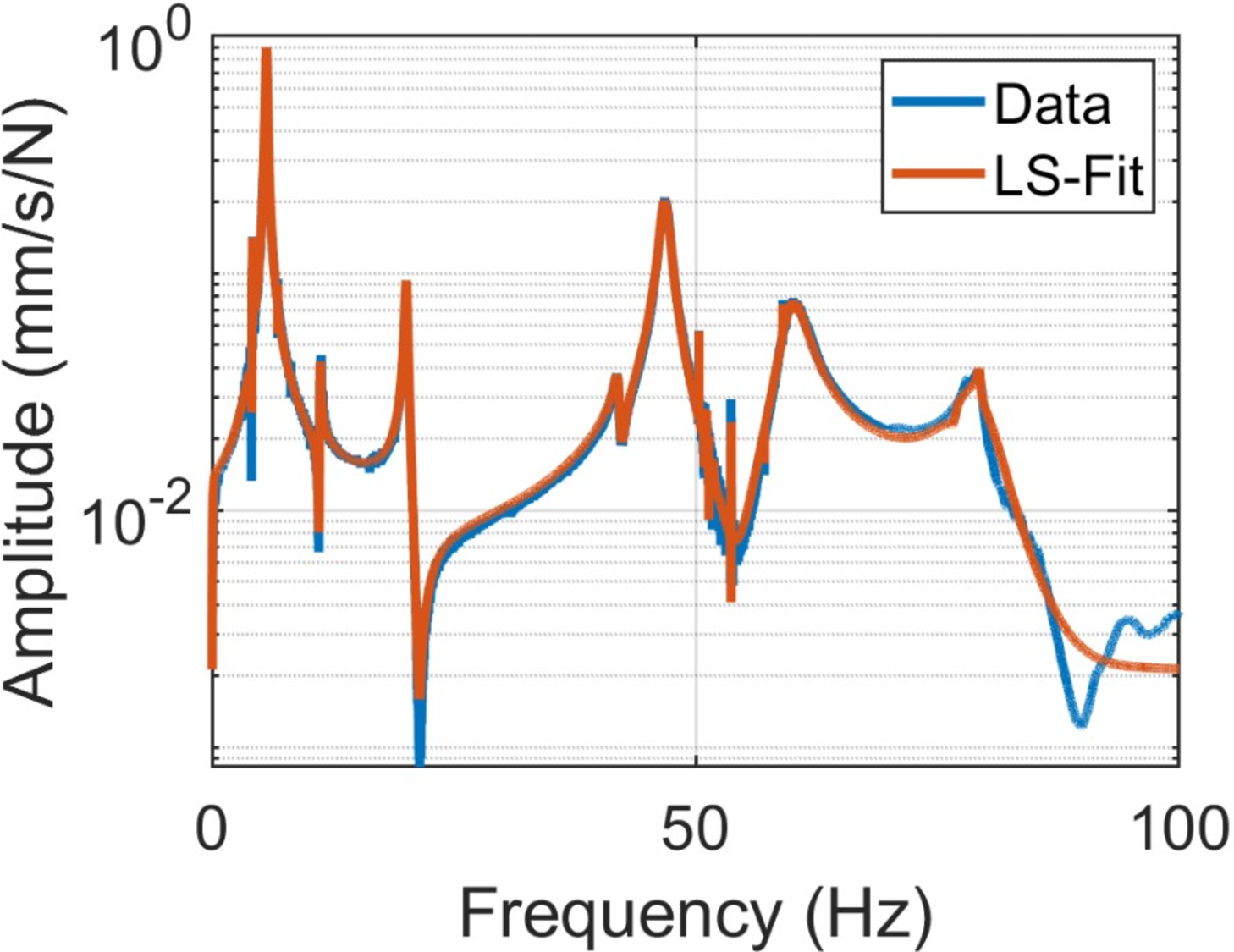}} &
       
     \subfloat[Load Case 6 : 0.0016 N $\mathrm{RMS}$]{\includegraphics[width=0.38\textwidth]{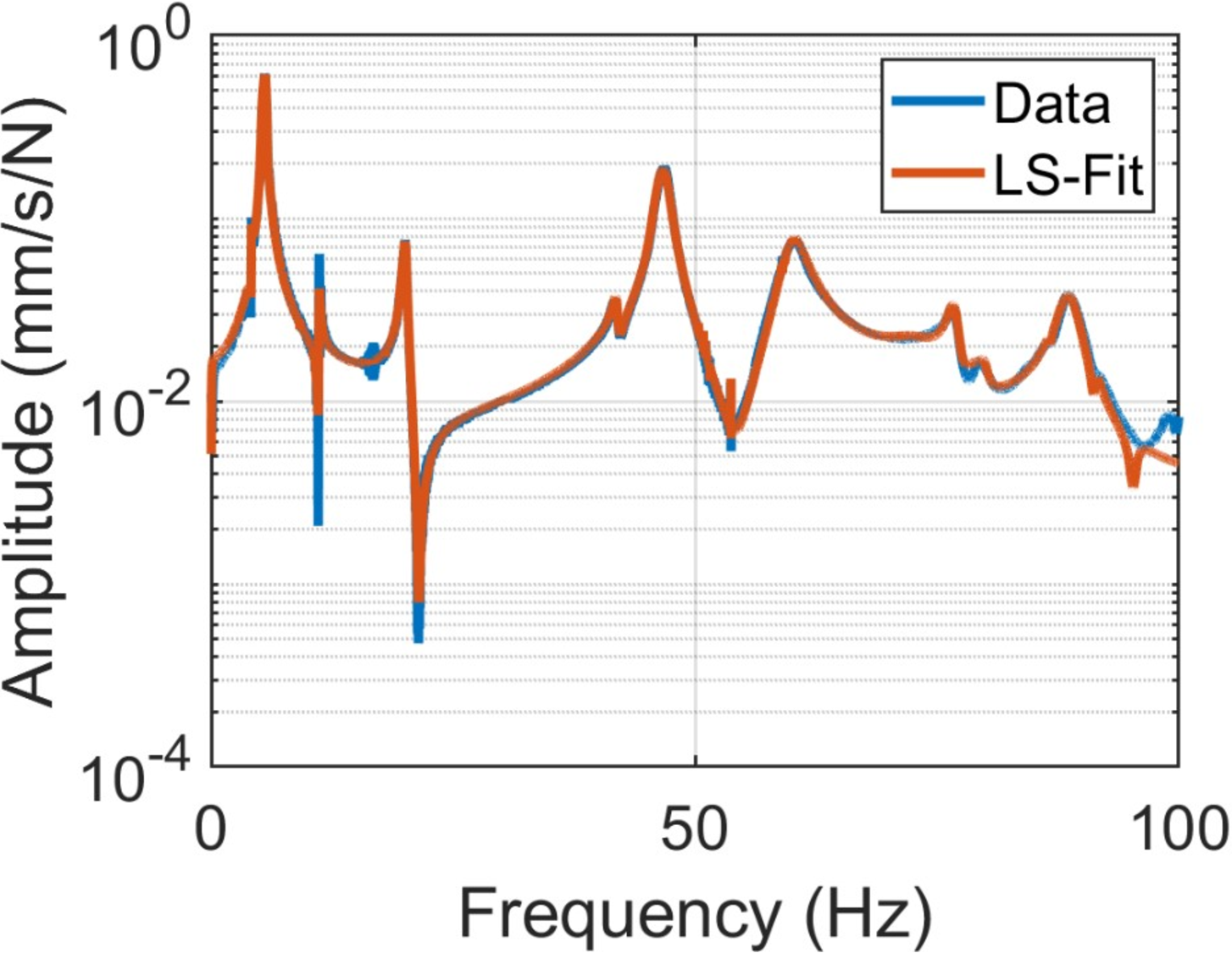}}\\

    \subfloat[Load Case 10 : 0.006 N  $\mathrm{RMS}$]{\includegraphics[width=0.38\textwidth]{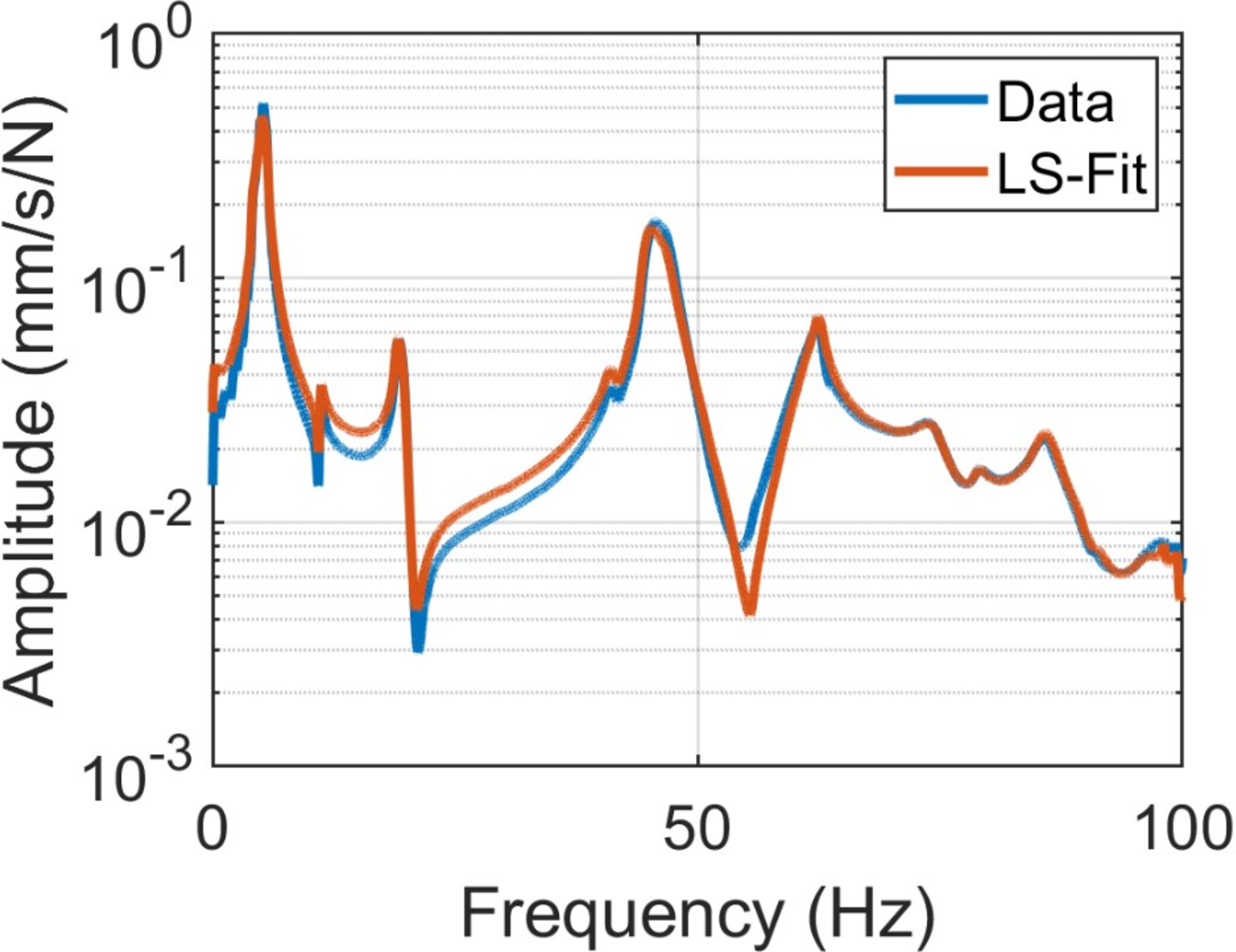}}&

    \subfloat[Load Case 15 : 0.01 N  $\mathrm{RMS}$]{\includegraphics[width=0.38\textwidth]{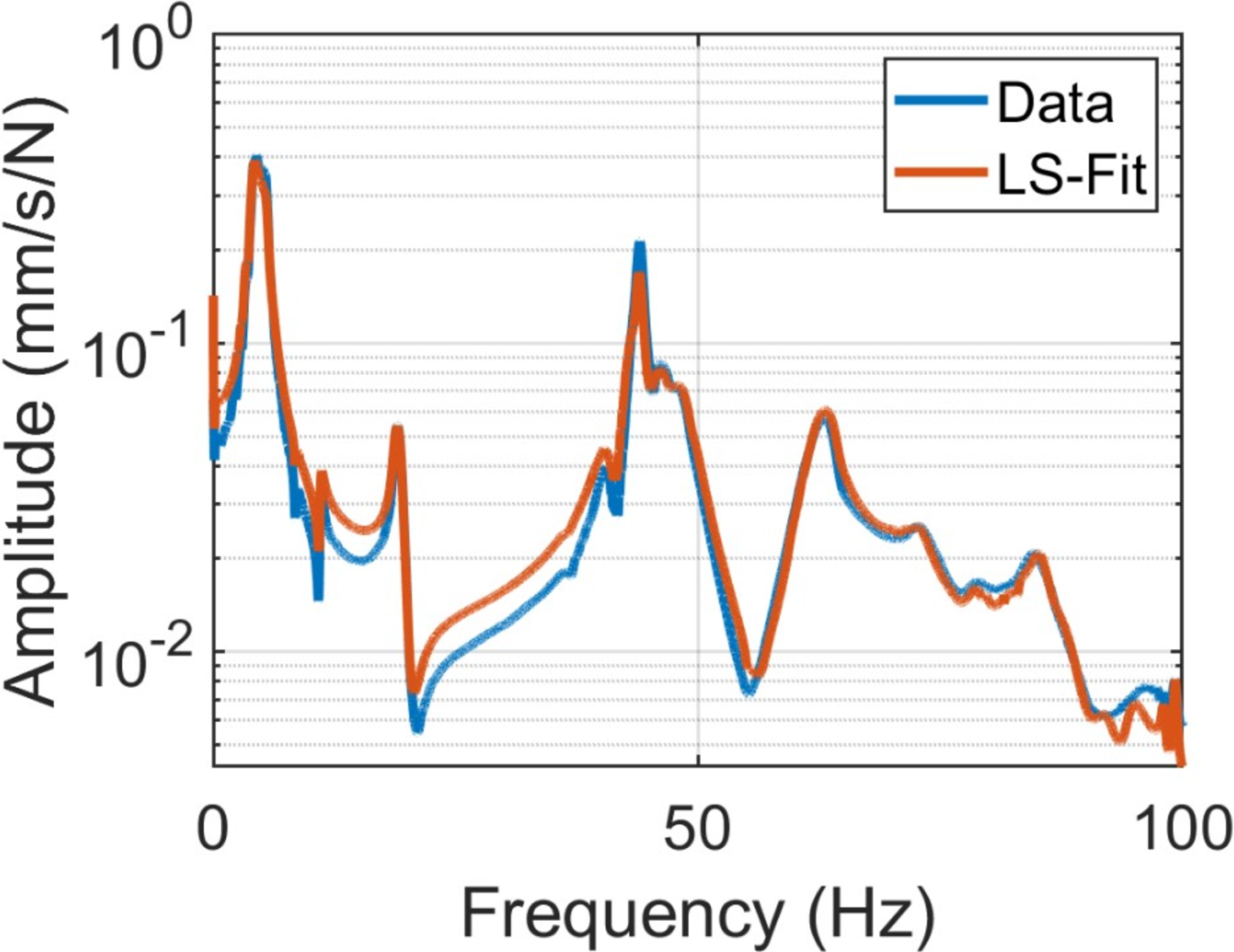}}\\
    
\end{tabular}
\caption{Least-Squares Fit (LS-Fit) using Vector Fitting with reduced order r = 200}
\label{NP_TFs}
\end{figure}

The objective in the non-parametric analysis is to construct a strictly proper degree-$r$ rational approximant $\tilde{H}(s)$ defined as the ratio of a numerator polynomial $n(s)$ and a denominator polynomial $d(s)$:
\begin{equation}\label{eq:6}
\tilde{H}(s) = \frac{n(s)}{d(s)} = \frac{\alpha_0 + \alpha_1 s + ... + \alpha_{r-1}s^{r-1}}{\beta_0 + \beta_1 s + ... + \beta_{r-1}s^{r-1} + s^r},
\end{equation}
where the unknown coefficients $\alpha_i$ and $\beta_i$ are chosen to minimize a weighted least squares (LS) error, $\varepsilon_{LS}$, between the measured data and the approximant: In other words, choose $\alpha_i$ and $\beta_i$ so that  
\begin{equation}\label{eq:7}
\varepsilon_{LS} = \sum_{i=1}^{N_s} w_i \left| H(s_i) - \frac{n(s_i)}{d(s_i)} \right|^2 \rightarrow \min,
\end{equation}
where $w_i$ are the associated sample weights. This nonlinear LS problem is solved via a sequence of linear LS problems using the Sanathanan-Koerner (SK) iteration scheme \cite{sk1963}. At the $k^{\text{th}}$ step, the denominator from the previous step, $d^{(k)}(s)$, is used as a known weighting function to solve for the updated polynomials:
\begin{equation}\label{eq:8}
\sum_{i=1}^{N_s} w_i \left| \frac{d^{(k+1)}(s_i)H(s_i) - n^{(k+1)}(s_i)}{d^{(k)}(s_i)} \right|^2 \rightarrow \min.
\end{equation}
Now, the cost function in Eq.~\eqref{eq:8} is linear in the unknowns $d^{(k+1)}$ and $n^{(k+1)}$.
To ensure numerical stability during these iterations, Vector Fitting~\cite{Gustavsen1999} utilizes the barycentric form of the rational approximant rather than the polynomial form where the approximant is expressed as 
\begin{equation}\label{eq:9}
\tilde{H}^{(k)}(s) = \frac{\tilde{n}^{(k)}(s)}{\tilde{d}^{(k)}(s)} = \frac{\sum_{i=1}^r \frac{\phi_i^{(k)}}{s-\lambda_i^{(k)}}}{1 + \sum_{i=1}^r \frac{\psi_i^{(k)}}{s-\lambda_i^{(k)}}},
\end{equation}
where $\phi_i^{(k)}\in \mathbb C$ and $\psi_i^{(k)}\in \mathbb C$ are called the barycentric coefficients and $\lambda_i^{(k)}$ the barycentric nodes at the $k^{\text{th}}$ iteration. It is easy to show that $\tilde{H}^{(k)}(s)$ in Eq.~\eqref{eq:9} is a rational function of degree-$r$. 
The Vector Fitting algorithm updates the barycentric nodes at each step by choosing the zeros of the updated denominator $\tilde{d}^{(k)}(s)$ as the starting poles for the next iteration. Once the iterations converges, $\tilde{d}^{(k)}(s) \rightarrow 1$, and the final approximant  is obtained in the so-called pole residue form, i.e., $$\Tilde{H}(s) = \sum_{i=1}^r\frac{\phi_i}{s-\lambda_i},$$
where $\lambda_i$ are the final poles of the rational approximant and $\phi_i$ the corresponding residues.  We emphasize that one can equivalently write this pole-residue form in the usual state-space form, 
\begin{equation} \label{eq:10}
    \Tilde{H}(s) = \sum_{i=1}^r\frac{\phi_i}{s-\lambda_i}=\Tilde{{C}}(s{I}-\Tilde{{A}})^{-1}\Tilde{{B}}     \hspace{2mm}\text{where}   \hspace{2mm} 
    \Tilde{{A}} \in \mathbb{R}^{r \times r},  \Tilde{{B}} \in \mathbb{R}^r, \text { and } \Tilde{{C}} \in \mathbb{R}^{1 \times r}.
\end{equation}
The experimental FRF data were used to develop a rational approximant for the boom system as 
in Eq.~\eqref{eq:10}. In this case, a set of $N_s=12800$ FRF samples was fitted via Vector Fitting using unity weights ($w_i=1$) via the ``Fast, Relaxed Vector-Fitting (FRVF)" implementation in Matlab \cite{Gustavsen2006,Deschrijver2008} to obtain a rational approximant of the form in Eq.~\eqref{eq:10} with order $r=200$.  To further enhance the quality of the fit at higher frequencies, Gaussian smoothing was applied to the data prior to Vector Fitting. This was done for all $15$ load-case models. Load cases $1$, $6$, $10$ and $15$ shown in Fig. \ref{NP_TFs} demonstrate that the models are an excellent fit for the actual experimental data. The models successfully match all the  key peak locations of the system across the tested frequency range.

Experimental validation of the derived rational approximation models was conducted using various excitation profiles such as sinusoidal, triangular and square waveforms. Force reconstruction was performed via frequency-domain inversion, where the estimated input force, $F_{\mathrm{test}}(s)$, is derived from the measured tip velocity, $V_{\mathrm{test}}(s)$, and the system transfer function according to  the frequency-domain equality 
$ V_{\mathrm{test}}(\imath \omega) = \tilde{F}_{\mathrm{test}}(\imath \omega)\tilde{H}_{\mathrm{ref}}(\imath \omega)$. Here, $\tilde{H}_{\mathrm{ref}}(\imath \omega)$ represents the transfer function obtained via Vector Fitting of the periodic chirp reference data. Then, the inverse Fourier transform of $\tilde{F}_{\mathrm{test}}(\imath \omega)$ yields the force signal, $\tilde{f}_\mathrm{test}(t)$ in the time domain. A comparison between estimated signal, $\tilde{f}_\mathrm{test}(t)$ and the measured signal, ${f}_\mathrm{test}(t)$ is made across all tested signals.

\begin{figure}[htbp]
\centering\includegraphics[width=0.9\textwidth]{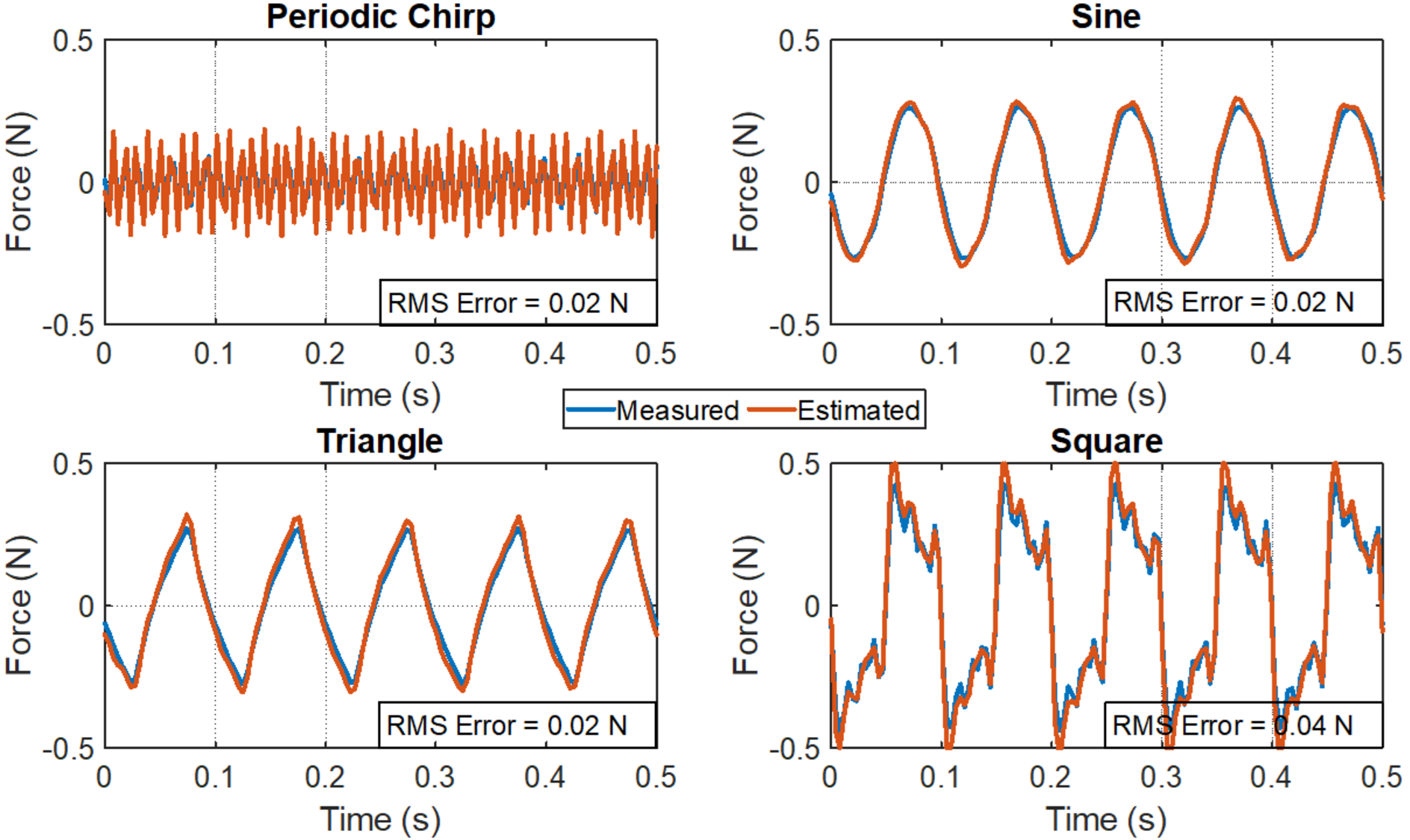}
\caption{Non-Parametric Modeling Validation Tests: Time Domain Input Estimation Error for Load Case 6 (0.0016 N RMS); RMS Errors between the Measured and Estimated signals inset in the subplots}
\label{NP_Freq_Test_Signals_time}
\end{figure}

The validation tests focused on Load Case 6 (0.0016 N RMS) have been shown in Fig. \ref{NP_Freq_Test_Signals_time}. Periodic signals were tested at a frequency of 10 Hz at the same RMS force level. The model exhibited strong estimation performance, particularly for the sine and triangle signals, where the estimated force closely matched the measured force. However, estimation errors increased for signals with sharp transitions, such as the square signal. RMS errors for each case have been shown as insets in the figures.

\begin{figure}[htbp]
\centering
\begin{tabular}{c c}

     \subfloat[Shift in first bending mode frequency with increase in load]{\includegraphics[width=0.45\textwidth]{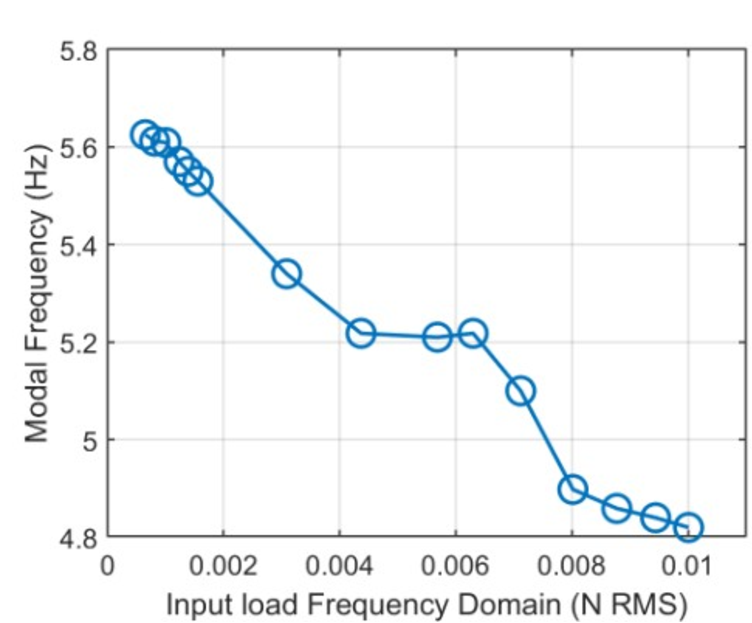}\label{modeshift}} &
       
     \subfloat[Cross-validation Error : Periodic Chirp Reference Signal All Combinations of Load Cases]{\includegraphics[width=0.5\textwidth]{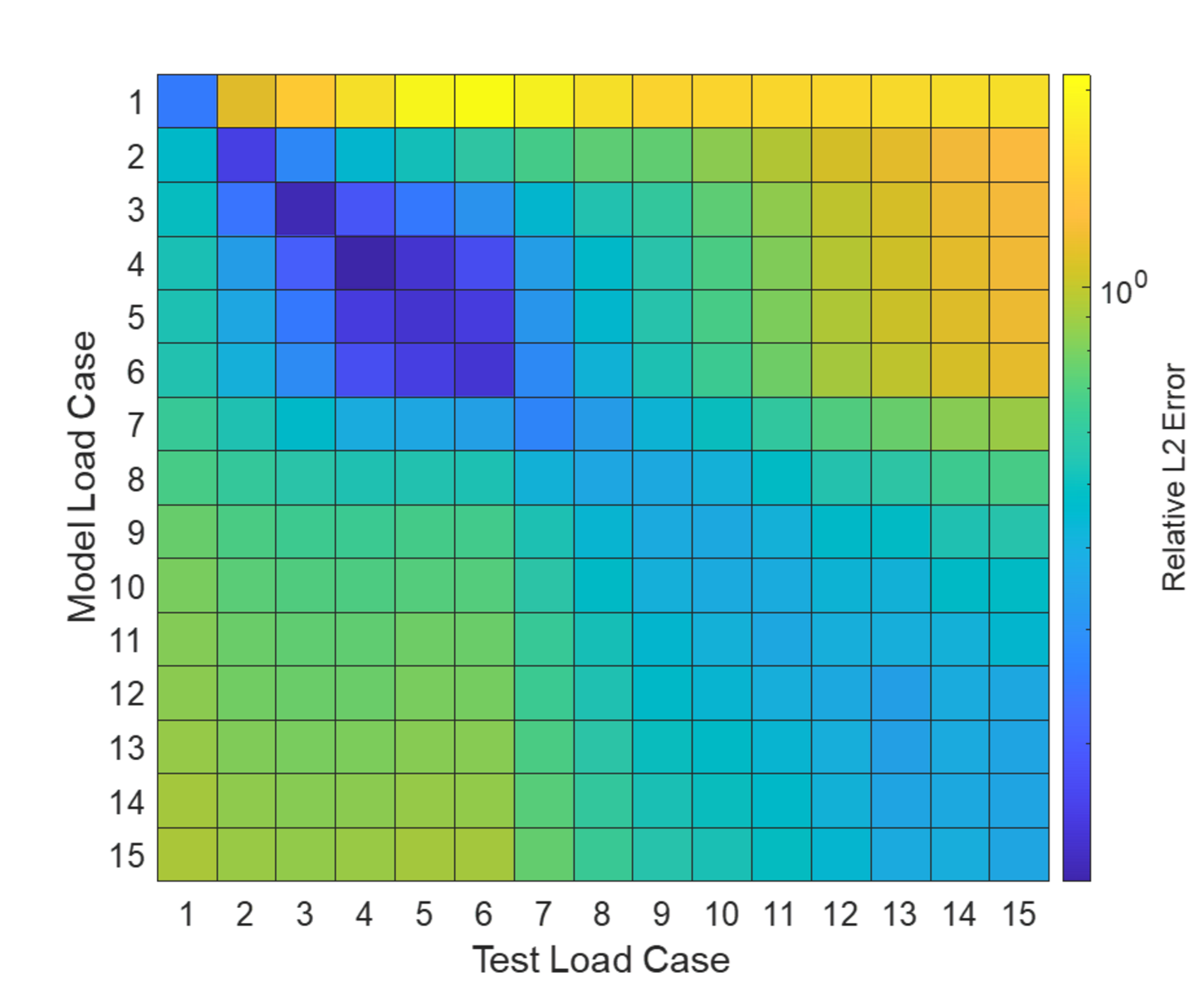}\label{cross_valid_all}}

\end{tabular}
\caption{Frequency Domain Inversion: Other Load Levels and the Need for Parametric Modeling}
\label{allres}
\end{figure}

As discussed earlier, there is a change in the dynamic behavior of the boom system with an increase in the input load amplitude. This load dependent behavior is observed in Fig.~\ref{modeshift} where the first bending mode frequency decreases from 5.6 to 4.8 Hz as the load increases from case 1 to 15. This shift in frequency is seen more clearly in Fig.~\ref{Noisy_to_Clean_TFs}.

We define a relative least-squares error ($\mathrm{E}^\mathrm{F}_\mathrm{L2}$) in the time domain to enable a scale-independent comparison of model fidelity across varying load intensities. For $N_t$ time-domain samples, 
\begin{equation}\label{eq:11}
\mathrm{E}^\mathrm{F}_\mathrm{L2}=\sqrt{\frac{\sum_{j=1}^{N_t}\left\|\tilde{f}_{\text{test}}(t_j)-f_{\text {test }}(t_j)\right\|_2^2}{\sum_{j=1}^{N_t}\left\|f_{\text{test }}(t_j)\right\|_2^2}},
\end{equation}
where $\tilde{f}_{\text{test}}$ and $f_{\text{test}}$ are the model predicted and measured test force-time histories respectively, at time step, $t_j$ . Fig. \ref{cross_valid_all} shows the Force error defined in Eq.~\eqref{eq:7} as a function of the test load case for all combinations of load cases (Model Load Case vs. Test Load Case), based on the Periodic Chirp reference signal. The large, green-yellow area above and below the diagonal signifies a large error, meaning that a non-parametric approach is not generalizable across the full range of operational loads. Therefore, a parametric model is necessary to properly capture the underlying system dynamics. 

\section{Parametric Modeling for Input Load Estimation}\label{parametric}

Fig. \ref{TF_all} showed a shift in the transfer functions with change in the input load.  Therefore, a non-parametric rational approximation approach will require retesting and refitting at each distinct load level due to the boom's highly load-dependent dynamic behavior, making it impractical for deployment across a continuous operational range. A parametric modeling framework directly addresses this limitation by constructing a single, unified transfer function model capable of capturing system behavior across all load levels simultaneously. Here, we generate a parametric model using the p-AAA algorithm \cite{Rodriguez2023} which extends the AAA algorithm \cite{Nakatsukasa2018} to multivariate problems. For this study, the dynamic response is parameterized by two distinct variables, frequency and input load intensity. We start with discussing the p-AAA methodology and apply it to our transfer-function dataset. Then, we illustrate its performance on various validation tests. Finally, a comparison is made with the non-parametric studies to evaluate its capability in predicting the load-dependent behavior of the boom.

\subsection{Input Force-based parameterization using p-AAA}

In contrast to the single variable case in Eq.~\eqref{eq:6}, we now consider rational approximation of a multivariate transfer function $H(s, p)$ that varies with both frequency ($s$), and the parameter (the load level) ($p$), defined earlier as the RMS value of the input force spectra in the $3  - 70 \text{ Hz}$ frequency range. We will collect the experimental data 
\begin{equation}\label{eq:12}
h_{j k}=H\left(s_j, p_k\right) \in \mathbb{C} \quad \text { for } \quad j=1, \ldots, N_s \quad \text { and } \quad k=1, \ldots, N_p ,
\end{equation}
where $N_s$ and $N_p$ are the total number of frequency samples and parameter values respectively. Similar to the single-variable case in non-parametric studies (in Eq.~\eqref{eq:9}), the rational approximant $\Tilde  H(s,p)$ is expressed in its two-variable barycentric form~\cite{ionita2014data,antoulas2012two}, namely 
 
\begin{equation}\label{eq:13}
    \Tilde{H}(s,p)=\frac{n(s,p)}{d(s,p)}=\left. \sum_{j=1}^l\sum_{k=1}^q\frac{\alpha_{jk}h_{jk}}{(s-\sigma_j)(p-\pi_k)} \middle/ \sum_{j=1}^l\sum_{k=1}^q\frac{\alpha_{jk}}{(s-\sigma_j)(p-\pi_k)} \right.,
\end{equation}
where $ \sigma_j$  and $\pi_k$  are to-be-determined points, subsets of the sets $\{s_j\}$  and $\{p_k\}$, respectively. The integers $l$ and $q$ represent the complexity (the order) of the approximant in the variables $s$ and $p$, respectively.

By construction, as long as $\alpha_{jk}\neq 0$, this barycentric form automatically interpolates the data at the selected support points $(\sigma_j, \pi_k)$ as the numerator weights are set to $\alpha_{jk}h_{jk}$. In other words,  $\Tilde{H}(\sigma_j,\pi_k) = {H}(\sigma_j,\pi_k)$ for 
$j=1,\ldots,l$ and $k=1,\ldots,q$.
So, in~\eqref{eq:13}, $\alpha_{jk}$ are the free parameters to be determined. Consequently, constructing the p-AAA approximant boils down to making two primary decisions: (1) choosing where to interpolate the data, i.e., selecting the support points $\sigma_j$ and $\pi_k$; and (2) determining how to compute the free parameters $\alpha_{jk}$ to best fit the remaining data.

To makes these decisions, the p-AAA algorithm relies on an iterative procedure that begins by partitioning the available data at each iteration $k$ into two disjoint, interpolation and least-squares fit, sets for both frequency and parameter samples: $\{\sigma\}$ and $\{\hat{\sigma}\}$ for frequency, and $\{\pi\}$ and $\{\hat{\pi}\}$ for the parameter: 
\begin{equation}\label{eq:14}
\begin{gathered}
\left[s_1, \ldots, s_{N_s}\right]=\overbrace{\left[\sigma_1, \ldots, \sigma_l\right]}^\text{Interpolation} \cup\overbrace{\left[\hat{\sigma}_1, \ldots, \hat{\sigma}_{{N_s}-l}\right]}^\text{Least-Squares Fitting} \\
\left[p_1, \ldots, p_{N_p}\right]=\left[\pi_1, \ldots, \pi_q\right] \cup\left[\hat{\pi}_1, \ldots, \hat{\pi}_{{N_p}-q}\right]
\end{gathered}
\end{equation}

 In each iteration, the p-AAA algorithm employs a greedy selection method to choose the next data pair, $(\sigma_l, \pi_q)$, to add to the interpolation set by finding the pair with the maximum current error:
\begin{equation}\label{eq:15}
\left(\sigma_l, \pi_q\right)=\argmax_{\left(s_j, p_k\right)}\left|h_{j k}-\tilde{H}\left(s_j, p_k\right)\right|.
\end{equation}
 As the algorithm proceeds through iterations, the interpolation set grows larger while the remaining set of data available for least-squares fitting shrinks. The coefficients, $\alpha_{jk}$, are then determined by minimizing a linearized least-squares (LS) error in the remaining data. This leads to a linearized LS problem solved in every step:
\begin{equation}\label{eq:16}
\tilde{H}=\underset{\hat{H}=n / d}{\arg \min } \sum_{j, k}\left|H\left(s_j, p_k\right) d\left(s_j, p_k\right)-n\left(s_j, p_k\right)\right|^2.
\end{equation}
The p-AAA iterations are run until the maximum error in~\eqref{eq:15} falls below a prespecified tolerance. For details on p-AAA and its recent extensions, we refer the reader to~\cite{Rodriguez2023,balicki2025multivariate,balicki2025scattered}. For other parametric data-driven modeling approaches that also leverage the barycentric form, we refer the reader to
\cite{antoulas2012two,antoulas2025loewner,Ionita2014}.

\begin{figure}[h!]
\centering
\includegraphics[width=0.8\textwidth]{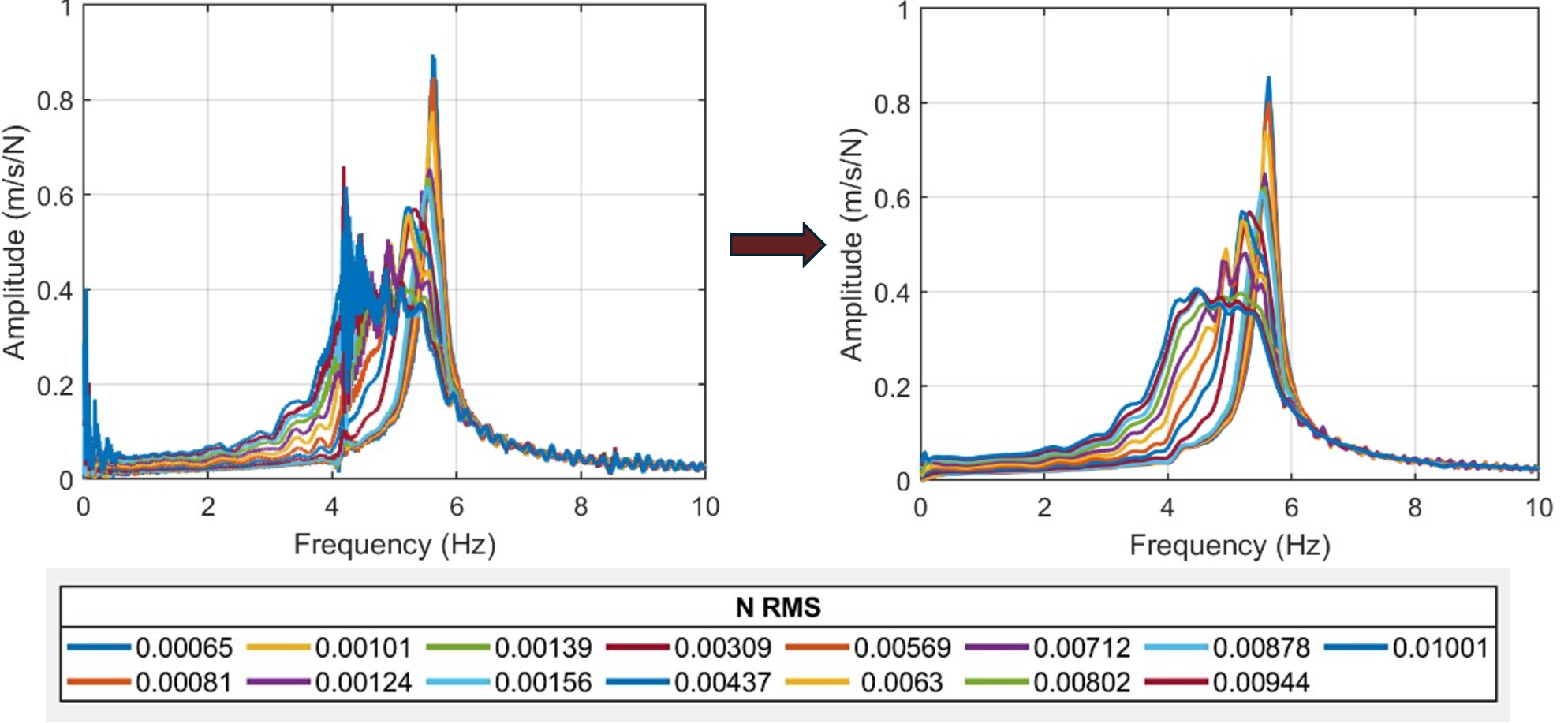}
\caption{Smoothing and Down-sampling for better Transfer-Function Fitting}
\label{Noisy_to_Clean_TFs}
\end{figure}

The combined experimental transfer function data for all the load cases is dense and noisy. Therefore,  for an effective implementation of the p-AAA algorithm, smoothing and down-sampling of the raw transfer function data are applied (shown in Fig. \ref{Noisy_to_Clean_TFs}) to reduce the computational time and improve model quality. Here, we use Gaussian smoothing because it selectively preserves the overall shape of the peaks in the transfer function, ensuring the resulting parametric model accurately captures the system's resonant behavior despite the decimation. The transfer-function data was down-sampled by a factor of 8 to $N_s=1600$ frequency samples. 

The p-AAA algorithm automatically selects the complexity (orders $l$ and $q$) of the approximant and the interpolation points, and runs iteratively until a pre-specified error tolerance is reached. p-AAA increases the order of the model with every iteration. While we originally targeted an order $l = 200$ (in frequency) p-AAA model to have a fair comparison with the non-parametric vector fitting model, the algorithm found the fit for the given tolerance at order $l$ = 198. Stopping the algorithm at this precise iteration ensures an optimally sized model and avoids overfitting to measurement noise.

\begin{figure}[ht]
\centering
\begin{tabular}{c c c}

     \subfloat[Load Case 1]{\includegraphics[width=0.38\textwidth]{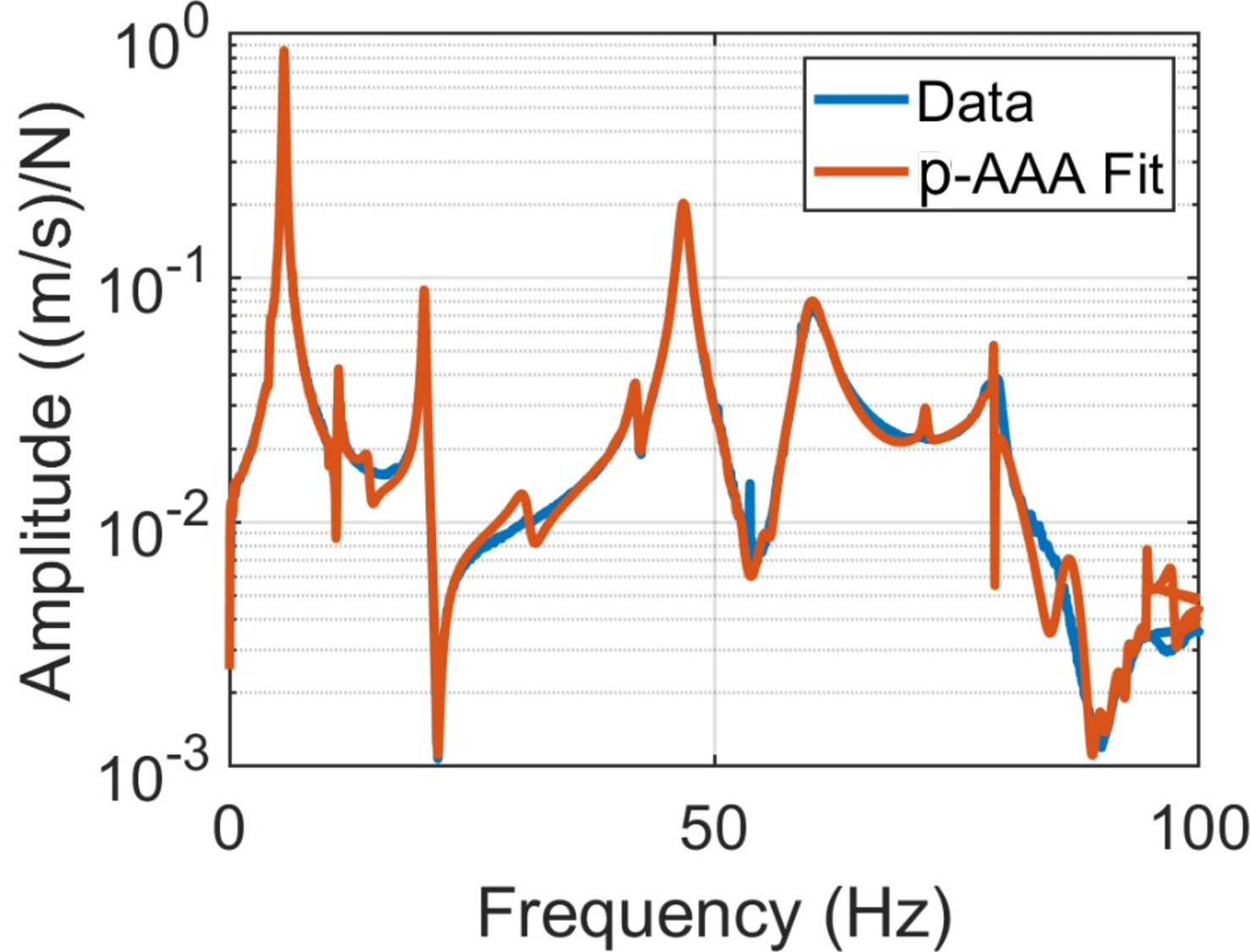}} &
       
     \subfloat[Load Case 6]{\includegraphics[width=0.38\textwidth]{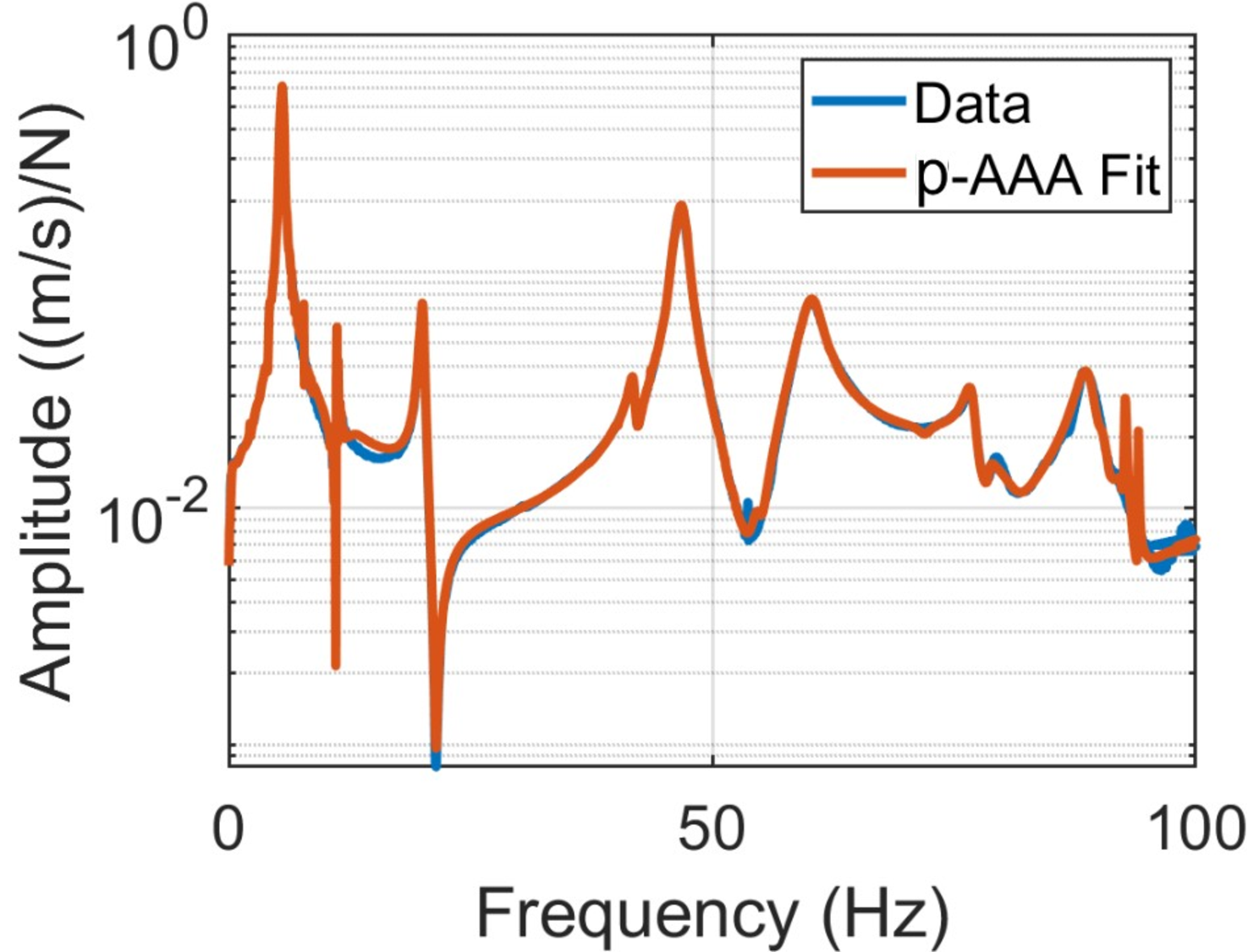}}\\

    \subfloat[Load Case 10]{\includegraphics[width=0.38\textwidth]{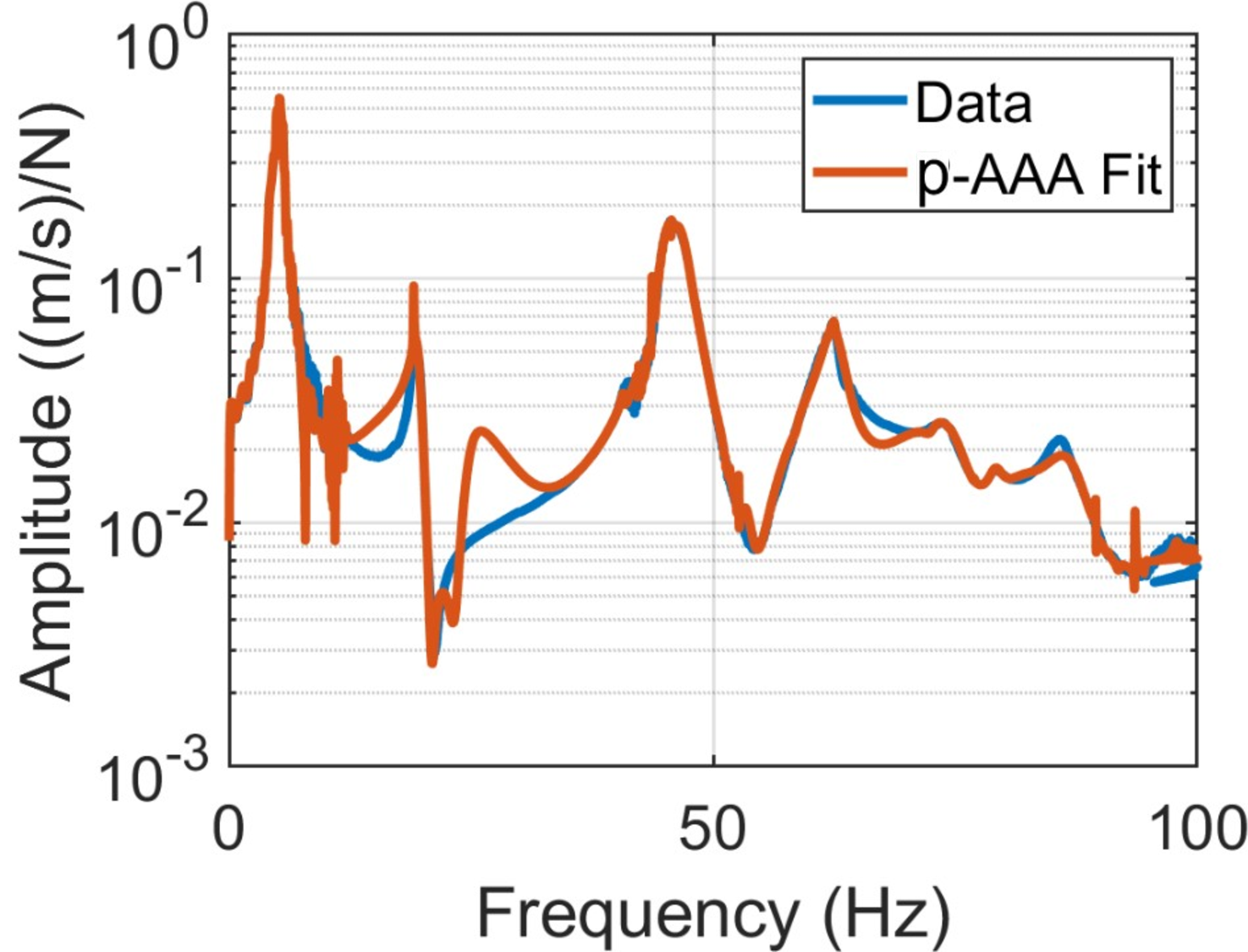}}&

    \subfloat[Load Case 15]{\includegraphics[width=0.38\textwidth]{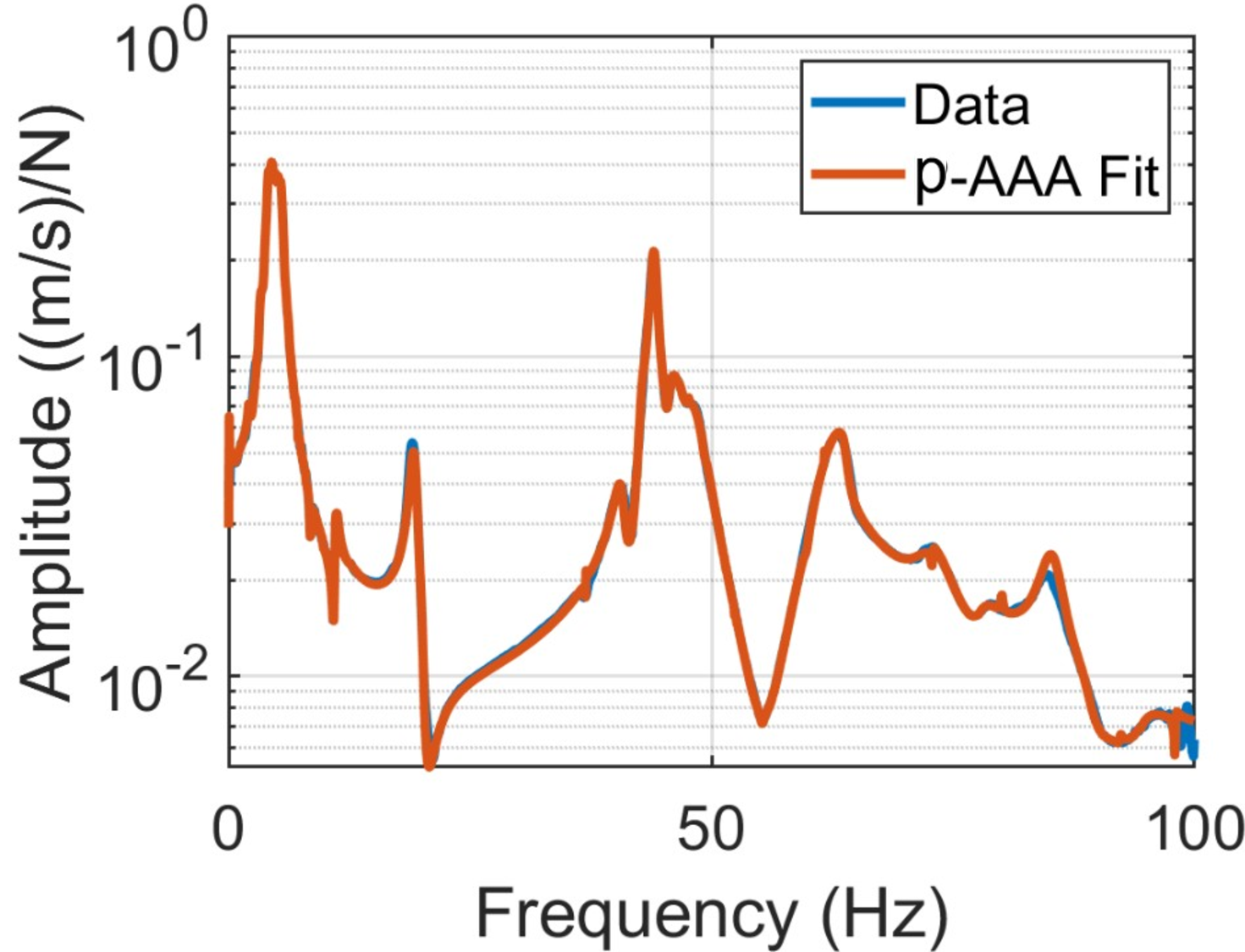}}\\

\end{tabular}
\caption{Parametric Modeling: Transfer Function Fit using p-AAA with order $l$=198 in frequency ($s$) and order $q$=14 in parameter index ($p$)}
\label{P_TFs}
\end{figure}

Thus, in summary, the p-AAA algorithm was used to construct a parametric rational approximation, $\Tilde{H}(s,p)$ (as shown in Eq.~\eqref{eq:9}), for this smoothened and reduced set of data (with $N_s=1600$ frequency samples and $N_p=15$ parameter samples). Then, for a given test value of the parameter $\hat{p}$, substituting it into Eq.~\eqref{eq:9}) yields a standard frequency response function, $\Tilde{H}(s,\hat{p})$, meaning that the model can be tested for different parameter values without additional testing for the new parameter. The results for load cases 1, 6, 10 and 15 are shown in Fig. \ref{P_TFs}. The plots demonstrate an excellent overall agreement between the experimental data and the p-AAA parametric model fit across the different load levels. The p-AAA fit successfully captures both the amplitude and the frequency of the resonant peaks across all four load cases. Finally, they show that a single, unified p-AAA model is able to adapt and accurately reflect the load-dependent changes without requiring separate models for each test case.

\subsection{Validation Tests}
\vspace{3 mm}

The derived 
multivariate rational approximation model is validated by applying various test signals such as sine, triangle and square, to the boom and then using parametric (p-AAA) models for frequency domain inversion to estimate the applied force against the directly measured force. Similar to the non-parametric (single variable) case shown earlier, the input force, $\tilde{F}_{\mathrm{test}}(\imath \omega)$, is estimated from the measured tip velocity, $V_{\mathrm{test}}(\imath \omega)$, and the parametric transfer function obtained from p-AAA fitting of the periodic chirp reference data, $\tilde{H}_{\mathrm{ref}}(\imath \omega,\hat{p})$ evaluated at the test parameter $\hat{p}$. This is done
as in Section~\ref{finv} using the frequency-domain equality 
$ V_{\mathrm{test}}(\imath \omega) = \tilde{F}_{\mathrm{test}}(\imath \omega)\tilde{H}_{\mathrm{ref}}(\imath \omega,\hat{p})$. Then, the inverse Fourier transform of $\tilde{F}_{\mathrm{test}}(\imath \omega)$ yields the force signal, $\tilde{f}_\mathrm{test}(t)$ in the time domain. 

\begin{figure}[htbp]
\centering
\includegraphics[width=0.9\textwidth]{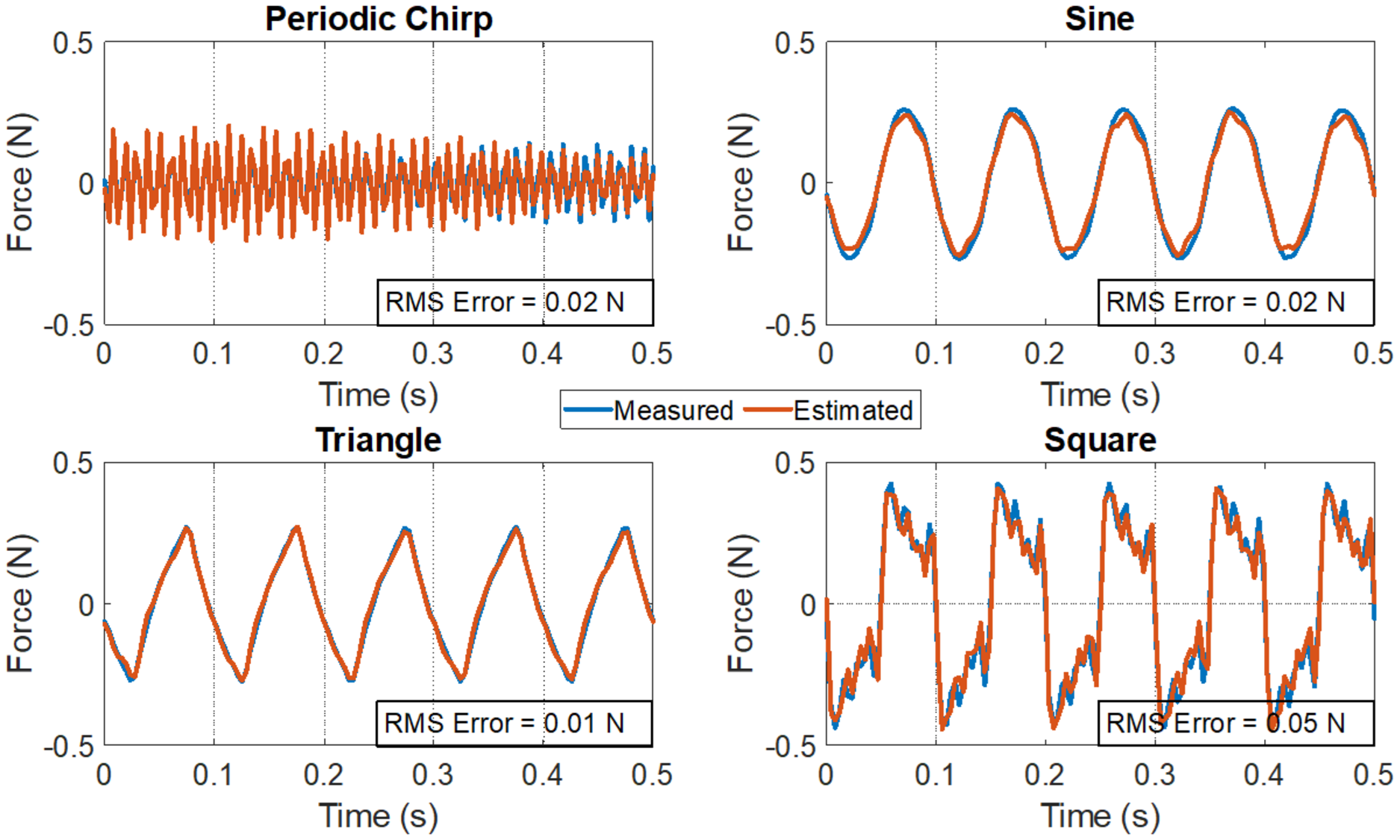}
\caption{Parametric Modeling Validation Tests: Time Domain Input Estimation Error for Load Case 6 (0.0016 N RMS); RMS Errors between the Measured and Estimated signals inset in the subplots}
\label{P_Freq_Test_Signals}
\end{figure}

A comparison between estimated signal, $\tilde{f}_\mathrm{test}(t)$ and the measured signal, ${f}_\mathrm{test}(t)$ is made across all tested signals. The parametric frequency domain model was observed to provide accurate force estimation over the full test range and the results for Load Case 6 ($p$ = 0.0016 N RMS) are shown in Fig. \ref{P_Freq_Test_Signals}. The estimated force accurately tracks both the amplitude and phase of the periodic chirp, sine, and triangle waveforms with small RMS errors ranging from 0.01 N to 0.02 N. For the Square wave, the model successfully captures the primary macroscopic force signal but exhibits a slightly higher RMS error of 0.05 N.

\subsection{Comparison with non-parametric analysis}
\vspace{3 mm}

\begin{figure}[h!]
\centering
\includegraphics[width=0.9\textwidth]{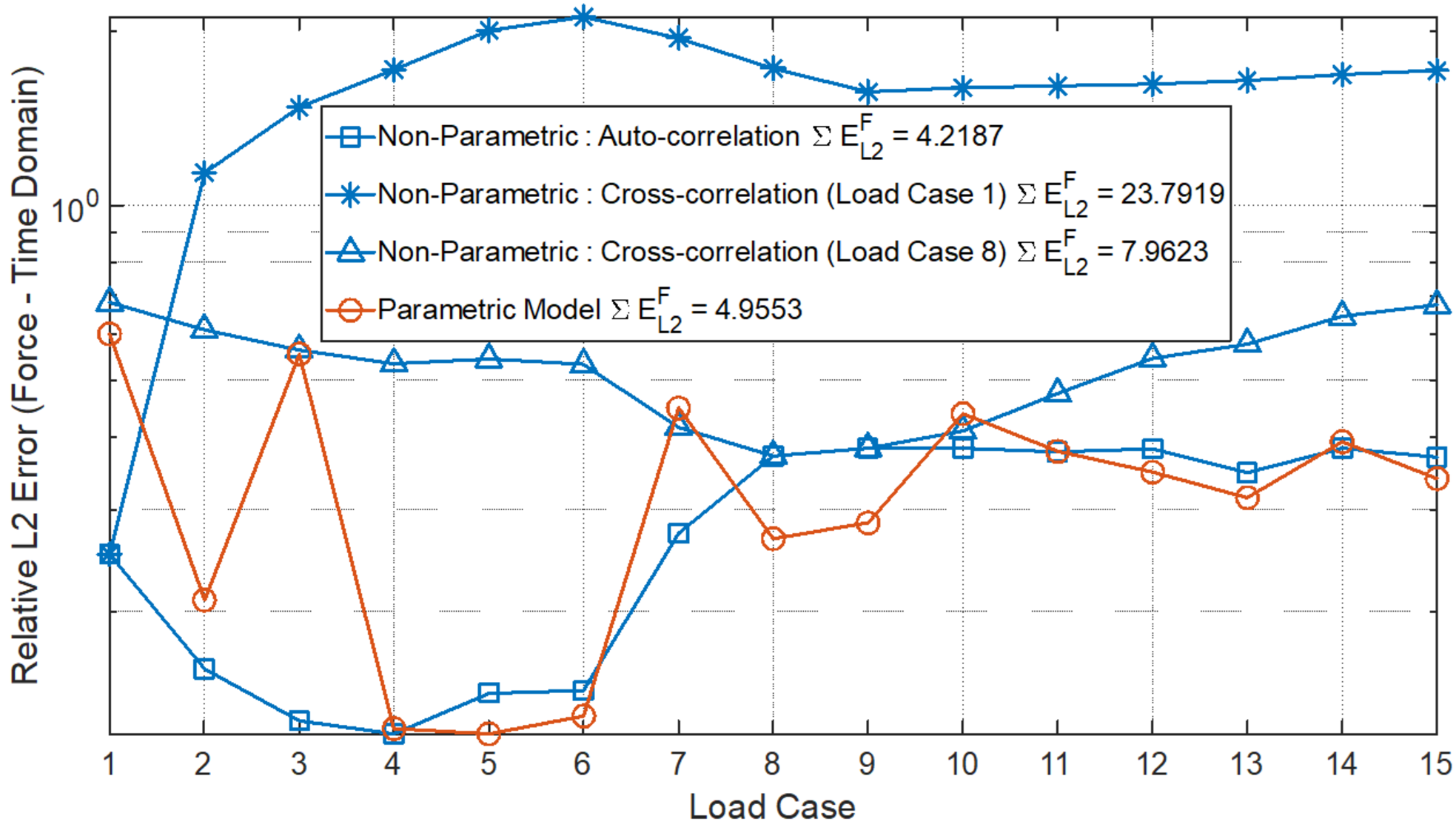}
\caption{Non-Parametric vs Parametric Analysis : Time Domain Reference Signal (Periodic Chirp) Estimation Error}
\label{PvsNP_PC}
\end{figure}

Fig. \ref{PvsNP_PC} shows a comparison of model performance across varying load levels for the 
periodic chirp reference signal. First, evaluating a non-parametric model trained on Load Case 1 (0.0006 N RMS) across all experimental load levels reveals an increase in error as the test load increases from the model load level. This significant decrease in performance demonstrates that a model fitted at a single, low excitation level cannot accurately generalize to higher load levels. Furthermore, the cross-validation results for a non-parametric model trained on Load Case 8 (0.004 N RMS) showed lowest total error across the entire test load range. However, despite its improved overall performance compared to the non-parametric model trained on the Load Case 1, it exhibits the same fundamental limitation where the error is minimized only at its specific training load level and increases when extrapolated to other excitation levels. The auto-correlation for the non-parametric analysis achieved a lower Force RMS Error compared to the cross-relation results across all load cases. This outcome is expected because auto-correlation evaluates the model against the exact same test data used to train it, resulting in an optimal fit for that single load point. However, this is impractical for real-world applications as relying on this would require re-testing at every distinct load parameter encountered during operation. 

On the other hand, the parametric model maintains a consistently low error between the measured and estimated input forces across all 15 load cases. In addition to this, the total error across all load levels for the parametric model ($\mathrm{E_{L2}^F}$ = 4.96) is much lower compared to the total error for the best non-parametric model based on Load Case 8 ($\mathrm{E_{L2}^F}$ = 7.96), thereby confirming that parametric modeling performs better than the non-parametric cross-correlation models, for generalizing across the full range of input loads. 

\begin{figure}[ht]
\centering
\includegraphics[width=0.9\textwidth]{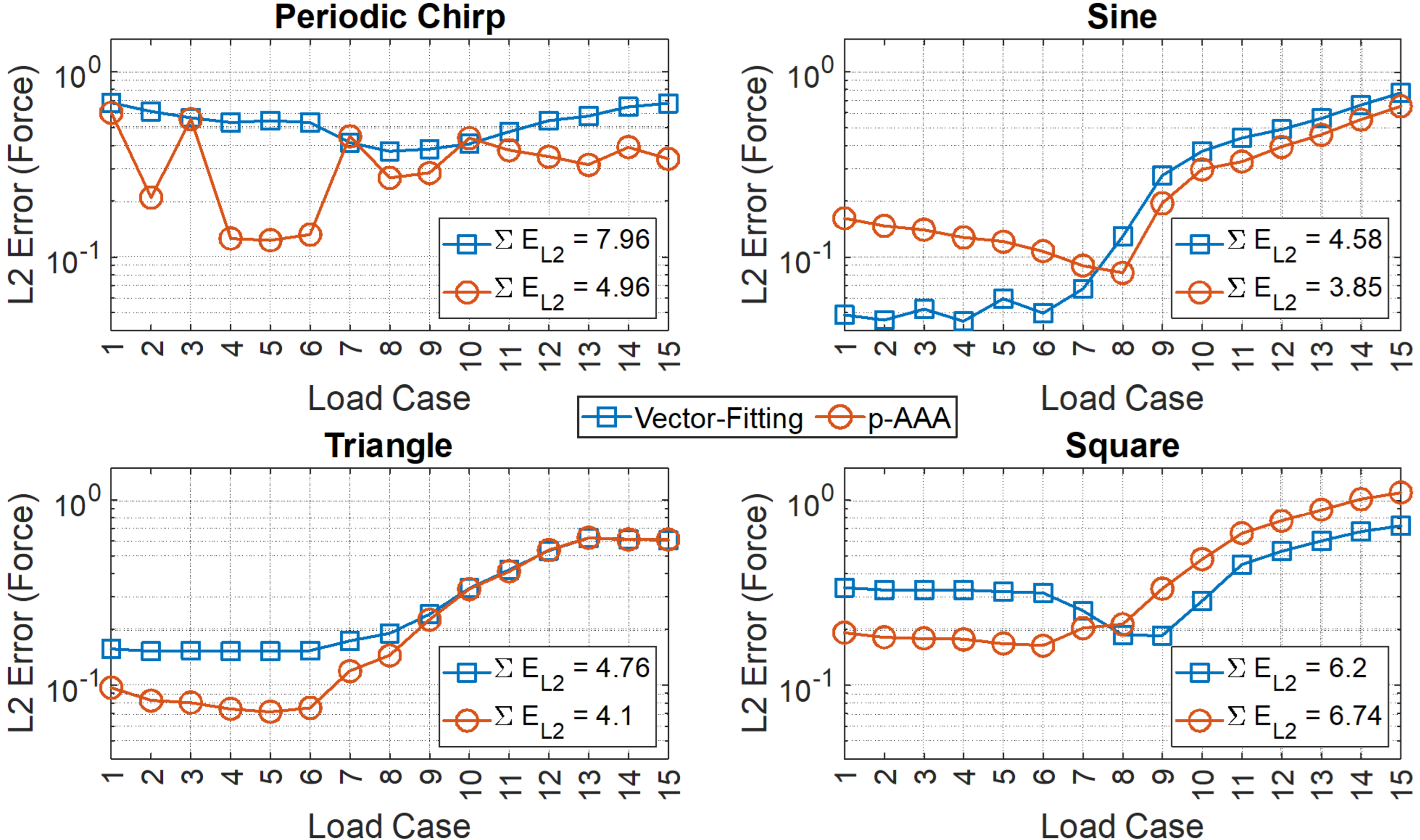}
\caption{Non-Parametric vs Parametric Analysis : Relative L2 Error for the force input (\textbf{$\mathrm{E}_\mathrm{L2}^\mathrm{F}$}) for the validation tests; Total Relative L2 Errors ($\Sigma \mathrm{E_{L2}^F}$) inset in the plots for model comparison, where the specific model load cases used for the Vector-Fitting results are 8, 2, 7, and 10 for the periodic chirp, sine, triangle, and square signals respectively}
\label{NP_vs_P_All_Test_Signals}
\end{figure}

A good criterion for evaluating the parametric framework is its ability to match or exceed the accuracy of load-specific models across diverse loading conditions. Fig. \ref{NP_vs_P_All_Test_Signals} presents this comparison by pitting the single parametric model against the most accurate load-specific non-parametric models for each test signal. Force signals were estimated for all combinations of model and test load cases for the non-parametric models and the total relative error ($\Sigma \mathrm{E_{L2}^F}$) was calculated for all the model load cases. Load cases 8, 2, 7 and 10 had the lowest total errors for the periodic chirp, sine, triangle and square signals respectively. As mentioned earlier, the total relative error between the estimated and measured reference input signals (periodic chirp) went down from 7.96 for the best non-parametric model to 4.96 for the parametric model, a reduction of 37.68\%. By quantifying the modeling error, this research successfully demonstrates that a single, unified parametric model can accurately capture the complex dynamics across the system's entire operational range. These results were supported by cross-validation for different test signals where the relative L2 error went down by 15.93\% from 4.58 to 3.85 and 13.86\% from 4.76 to 4.1 for the sine and triangle signals. The square signal was an exception where the total relative L2 error increased slightly by 8.7\% from 6.2 to 6.74 for the parametric model but this might be attributed to the sudden transitions inherent to the signal. Overall, these findings show that the parametric approach provides a generalizable force estimation solution across a continuous range of input loads.

\vspace{1 cm}

\section{Conclusion and Future Work}\label{conclusion}

This study demonstrates the need for parameterization when analyzing the structural dynamics of deployable booms. Non-parametric modeling provided good results at discrete load levels but failed to predict the system's evolving behavior with varying input loads. This problem was solved by using parametric frequency domain modeling via the p-AAA algorithm. The validation test results demonstrate that the p-AAA model yields a lower or comparable total relative error for the periodic chirp, sine, triangle and square signals compared to the non-parametric model. The methodology achieved a nearly 38\% reduction in total relative force estimation error for the reference signal. Having established this high-accuracy parametric model for the inverse problem, future work will focus on carrying out physical deployment tests  and directly utilizing this single p-AAA model to estimate the actual deployment loads. Additionally, the experimental framework will be refined to include precise input force matching and velocity-based parameterization to better capture dynamic variations and further improve the robustness of the load estimation. Ultimately, this parametric deployment load estimation framework will be utilized to predict the transient structural response of the satellite chassis during passive boom deployments in space. Its immediate application is the UPS-1 mission, where it will help assess mechanical shocks and protect sensitive onboard systems during orbital operations.

\section*{Acknowledgments}
We gratefully acknowledge the support from the Center for Space Science and Engineering Research at Virginia Tech (Space@VT) and VT Mechanical Engineering Department faculty Dr. Suyi Li for letting us use their lab space and equipment for the experimental work. We would also like to thank NASA Langley Research Center, Drs. Juan M. Fernandez and Matthew Chamberlain, for providing the boom samples for testing and flight. Gugercin's work was funded in part by US National Science Foundation grant DMS-2411141.

\bibliography{References}

\end{document}